\documentclass[11pt,twoside]{article}

\usepackage{amsfonts,epsfig,graphicx,subfigure}
\usepackage{amsmath,amssymb,amsthm}
\usepackage{fullpage}

\usepackage{fullpage}

\usepackage{epsf}
\usepackage{fancyheadings}
\usepackage{graphics,color}
\usepackage{epsf}
\usepackage{psfrag}

\usepackage{amsthm}
\usepackage{amsfonts}
\usepackage{amsmath}
\usepackage{amssymb}

\newlength{\widebarargwidth}
\newlength{\widebarargheight}
\newlength{\widebarargdepth}

\newcommand{\matsnorm}[2]{|\!|\!| #1 | \! | \!|_{{#2}}}

\newenvironment{carlist}
 {\begin{list}{$\bullet$}
 {\setlength{\topsep}{0in} \setlength{\partopsep}{0in}
  \setlength{\parsep}{0in} \setlength{\itemsep}{\parskip}
  \setlength{\leftmargin}{0.07in} \setlength{\rightmargin}{0.08in}
  \setlength{\listparindent}{0in} \setlength{\labelwidth}{0.08in}
  \setlength{\labelsep}{0.1in} \setlength{\itemindent}{0in}}}
 {\end{list}}

\newcommand{\bcar}{\begin{carlist}}
\newcommand{\ecar}{\end{carlist}}

\newcommand{\tracer}[2]{\ensuremath{\langle \!\langle {#1}, \; {#2}
\rangle \!\rangle}}

\newcommand{\op}[1]{\ensuremath{\operatorname{#1}}}

\makeatletter
\long\def\@makecaption#1#2{
        \vskip 0.8ex
        \setbox\@tempboxa\hbox{\small {\bf #1:} #2}
        \parindent 1.5em  
        \dimen0=\hsize
        \advance\dimen0 by -3em
        \ifdim \wd\@tempboxa >\dimen0
                \hbox to \hsize{
                        \parindent 0em
                        \hfil 
                        \parbox{\dimen0}{\def\baselinestretch{0.96}\small
                                {\bf #1.} #2
                                } 
                        \hfil}
        \else \hbox to \hsize{\hfil \box\@tempboxa \hfil}
        \fi
        }
\makeatother

\long\def\comment#1{}

\makeatletter
\def\@cite#1#2{[\if@tempswa #2 \fi #1]}
\makeatother

\makeatletter
\long\def\barenote#1{
    \insert\footins{\footnotesize
    \interlinepenalty\interfootnotelinepenalty 
    \splittopskip\footnotesep
    \splitmaxdepth \dp\strutbox \floatingpenalty \@MM
    \hsize\columnwidth \@parboxrestore
    {\rule{\z@}{\footnotesep}\ignorespaces
      #1\strut}}}
\makeatother

\newcommand{\bit}{\begin{itemize}}
\newcommand{\eit}{\end{itemize}}
\newcommand{\ben}{\begin{enumerate}}
\newcommand{\een}{\end{enumerate}}

\newcommand{\bear}{\begin{eqnarray}}
\newcommand{\eear}{\end{eqnarray}}

\newcommand{\inprod}[2]{\ensuremath{\langle #1 , \, #2 \rangle}}

\newcommand{\Exs}{\ensuremath{{\mathbb{E}}}}

\newcommand{\beq}{\begin{quotation}}
\newcommand{\enq}{\end{quotation}}

\newcommand{\estart}{\begin{equation}}
\newcommand{\eend}{\end{equation}}

\newcommand{\widgraph}[2]{\includegraphics[keepaspectratio,width=#1]{#2}}

\newcommand{\defn}{\ensuremath{:  =}}

\newcommand{\bec}{\begin{center}}
\newcommand{\enc}{\end{center}}

\newcommand{\beit}{\begin{itemize}}
\newcommand{\enit}{\end{itemize}}

\newcommand{\been}{\begin{enumerate}}
\newcommand{\enen}{\end{enumerate}}

\newcommand{\comsl}{\begin{slide}}
\newcommand{\comspor}{\begin{slide*}}

\newcommand{\comsld}[2]{\begin{slide}[#1,#2]}
\newcommand{\comspord}[2]{\begin{slide*}[#1,#2]}

\newcommand{\mendsl}{\end{slide}}
\newcommand{\mendspo}{\end{slide*}}

\newcommand{\real}{\ensuremath{{\mathbb{R}}}}

\DeclareMathOperator{\var}{var}

\DeclareMathOperator{\trace}{trace}

\theoremstyle{plain}

\newtheorem{theo}{Theorem}

\newtheorem{lem}{Lemma}[section]
\newtheorem{prop}{Proposition}[section]
\newtheorem{cor}{Corollary}

\theoremstyle{definition} 

\newtheorem{nota}{Notation}[section]
\newtheorem{de}{Definition}[section]
\newtheorem{exa}{Example}[section]
\newtheorem{as}{Assumption}[section]
\newtheorem{alg}{Algorithm}[section]

\newcommand{\btheo}{\begin{theo}}
\newcommand{\bde}{\begin{de}}
\newcommand{\ble}{\begin{lem}}
\newcommand{\bpr}{\begin{prop}}
\newcommand{\bno}{\begin{nota}}
\newcommand{\bex}{\begin{exa}}
\newcommand{\bcor}{\begin{cor}}
\newcommand{\spro}{\begin{proof}}
\newcommand{\bas}{\begin{as}}
\newcommand{\balg}{\begin{alg}}

\newcommand{\etheo}{\end{theo}}
\newcommand{\ede}{\end{de}}
\newcommand{\ele}{\end{lem}}
\newcommand{\epr}{\end{prop}}
\newcommand{\eno}{\end{nota}}
\newcommand{\eex}{\end{exa}}
\newcommand{\ecor}{\end{cor}}
\newcommand{\fpro}{\end{proof}}
\newcommand{\eas}{\end{as}}
\newcommand{\ealg}{\end{alg}}

\theoremstyle{plain}

\newtheorem{theos}{Theorem}
\newtheorem{props}{Proposition}
\newtheorem{lems}{Lemma}
\newtheorem{cors}{Corollary}

\theoremstyle{definition}
\newtheorem{exas}{Example}
\newtheorem{algs}{Algorithm}
\newtheorem{asss}{Assumption}
\newtheorem{defns}{Definition}

\newcommand{\btheos}{\begin{theos}}
\newcommand{\etheos}{\end{theos}}
\newcommand{\bprops}{\begin{props}}
\newcommand{\eprops}{\end{props}}
\newcommand{\bdes}{\begin{defns}}
\newcommand{\edes}{\end{defns}}
\newcommand{\blems}{\begin{lems}}
\newcommand{\elems}{\end{lems}}
\newcommand{\bcors}{\begin{cors}}
\newcommand{\ecors}{\end{cors}}
\newcommand{\bexs}{\begin{exas}}
\newcommand{\eexs}{\end{exas}}
\newcommand{\balgs}{\begin{algs}}
\newcommand{\ealgs}{\end{algs}}
\newcommand{\bass}{\begin{asss}}
\newcommand{\eass}{\end{asss}}

\def\real{\mathbb{R}}

\newcommand{\numobs}{\ensuremath{n}} 
\newcommand{\maxop}{\gamma} 

\newcommand{\Error}{\Delta} 
\newcommand{\nucnorm}[1]{\matsnorm{#1}{1}} 
\newcommand{\paramvar}{\Theta}
\newcommand{\og}{\ensuremath{\Theta^*}}
\newcommand{\thetastar}{\ensuremath{\og}}

\DeclareMathOperator{\rank}{rank}

\newcommand{\binprod}[2]{\ensuremath{\big \langle #1 , \, #2 \big \rangle}}

\newcommand{\pdim}{\ensuremath{p}}

\newcommand{\mprob}{\ensuremath{\mathbb{P}}}

\newcommand{\weakq}{R_q}

\newcommand{\ThetaHat}{\ensuremath{\widehat{\Theta}}}

\newcommand{\CovMat}{\ensuremath{\Sigma}}

\newcommand{\Sphere}[1]{\ensuremath{S^{#1-1}}}

\newcommand{\opnorm}[1]{\ensuremath{\matsnorm{#1}{\myop}}}
\newcommand{\frobnorm}[1]{\ensuremath{\matsnorm{#1}{F}}}

\newcommand{\sampcov}{\ensuremath{\frac{1}{\numobs} \SpecX^T \SpecX}}
\newcommand{\Sampcov}{\ensuremath{\Big(\sampcov\Big)}}

\newcommand{\ex}{E}

\newcommand{\weakball}{\ensuremath{\mathbb{B}_q(\weakq)}}

\newcommand{\ellzball}{\ensuremath{\mathbb{B}_0(R_0)}}

\newcommand{\mybigY}{\ensuremath{\vec{y}}}
\newcommand{\mybigX}{\ensuremath{X}}
\newcommand{\mybigW}{\ensuremath{\vec{\varepsilon}}}
\newcommand{\ThetaStar}{\ensuremath{\Theta^*}}

\newcommand{\mysmally}[1]{\ensuremath{y_{#1}}}
\newcommand{\mysmallX}[1]{\ensuremath{X_{#1}}}
\newcommand{\mysmallw}[1]{\ensuremath{\varepsilon_{#1}}}

\newcommand{\Xop}{\ensuremath{\mathfrak{X}}}
\newcommand{\Xopadj}{\ensuremath{\Xop^*}}
\newcommand{\lra}{\ensuremath{k}}
\newcommand{\lrb}{\ensuremath{\pdim}}

\newcommand{\Zsam}[1]{\ensuremath{Z_{#1}}}
\newcommand{\Noisesam}[1]{\ensuremath{W_{#1}}}

\newcommand{\bigN}{\ensuremath{N}}

\newcommand{\Xcov}[1]{\ensuremath{Z_{#1}}}
\newcommand{\Zout}[1]{\ensuremath{Y_{#1}}}
\newcommand{\Wout}[1]{\ensuremath{W_{#1}}}

\newcommand{\bigreg}{\ensuremath{\lambda_\bigN}}
\newcommand{\regpar}{\ensuremath{\lambda_n}}

\newcommand{\rdim}{\ensuremath{r}} 

\newcommand{\mata}{\ensuremath{U}} 
\newcommand{\matb}{\ensuremath{V}}

\newcommand{\mataspec}{\ensuremath{\widetilde{\mata}}}
\newcommand{\matbspec}{\ensuremath{\widetilde{\matb}}}

\newcommand{\SetOne}{\mathcal{M}}
\newcommand{\SetTwo}{\SetOne^\perp}

\newcommand{\row}{\ensuremath{\operatorname{row}}}
\newcommand{\col}{\ensuremath{\operatorname{col}}}

\newcommand{\PlError}{\ensuremath{\Delta}}

\newcommand{\myop}{\ensuremath{\operatorname{op}}}

\newcommand{\qpar}{\ensuremath{q}}
\newcommand{\radq}{\ensuremath{R_\qpar}}
\newcommand{\radz}{\ensuremath{R_0}}

\newcommand{\tmp}{m} 
\newcommand{\Kcomp}{\ensuremath{{K^c}}}
\newcommand{\mydefn}{\ensuremath{: \, =}}

\newcommand{\PlErrorA}{\ensuremath{\PlError'}}
\newcommand{\PlErrorB}{\ensuremath{\PlError''}}

\newcommand{\SpecX}{\ensuremath{X}}
\newcommand{\SpecY}{\ensuremath{Y}}
\newcommand{\SpecV}{\ensuremath{W}}
\newcommand{\sigmin}{\ensuremath{\sigma_{\operatorname{min}}}}
\newcommand{\sigmax}{\ensuremath{\sigma_{\operatorname{max}}}}

\newcommand{\Covera}{\ensuremath{\mathcal{A}}}
\newcommand{\covnuma}{\ensuremath{A}}
\newcommand{\Coverb}{\ensuremath{\mathcal{B}}}
\newcommand{\covnumb}{\ensuremath{B}} 

\newcommand{\mymin}{\ensuremath{m}}
\newcommand{\DelRotate}{\ensuremath{\Gamma}}

\newcommand{\KeySet}{\ensuremath{\mathcal{C}}}

\newcommand{\Tail}{\ensuremath{\mathcal{T}}}

\newcommand{\nvar}{\ensuremath{\nu}}

\newcommand{\arbrad}{\ensuremath{t}}
\newcommand{\SpecSet}{\ensuremath{\mathcal{R}}}
\newcommand{\conone}{\frac{1}{4}}
\newcommand{\contwo}{}

\newcommand{\matardim}{\ensuremath{\mata^\rdim}}
\newcommand{\matbrdim}{\ensuremath{\matb^\rdim}}

\newcommand{\Sah}{\ensuremath{T}} 

\newcommand{\BigMess}{\left(\sqrt{\frac{\lra}{\bigN}} +
\sqrt{\frac{\lrb}{\bigN}} \right)}

\begin{document}

  
  \begin{center}
    
    {\bf{\LARGE{Estimation of (near) low-rank
    matrices \\ with noise and high-dimensional scaling}}}
    
    \vspace*{.2in}

	    {\large{
		\begin{tabular}{ccc}
		  Sahand Negahban$^\star$ & & Martin
		  J. Wainwright$^{\dagger,\star}$ \\
		\end{tabular}
	    }}

	    \vspace*{.2in}

	    \begin{tabular}{c}
	      Department of Statistics$^\dagger$, and \\
	      Department of Electrical Engineering and Computer Sciences$^\star$ \\
	      UC Berkeley,  Berkeley, CA  94720
	    \end{tabular}

	    \vspace*{.2in}

	    \today

	    \vspace*{.2in}

	    \begin{tabular}{c}
	      Technical Report, \\
	      Department of Statistics,  UC Berkeley
	    \end{tabular}

\begin{abstract} 
  High-dimensional inference refers to problems of statistical
  estimation in which the ambient dimension of the data may be
  comparable to or possibly even larger than the sample size.  We
  study an instance of high-dimensional inference in which the goal is
  to estimate a matrix $\Theta^* \in \real^{k \times p}$ on the basis
  of $N$ noisy observations, and the unknown matrix $\Theta^*$ is
  assumed to be either exactly low rank, or ``near'' low-rank, meaning
  that it can be well-approximated by a matrix with low rank.  We
  consider an $M$-estimator based on regularization by the trace or
  nuclear norm over matrices, and analyze its performance under
  high-dimensional scaling.  We provide non-asymptotic bounds on the
  Frobenius norm error that hold for a general class of noisy
  observation models, and then illustrate their consequences for a
  number of specific matrix models, including low-rank multivariate or
  multi-task regression, system identification in vector
  autoregressive processes, and recovery of low-rank matrices from
  random projections.  Simulation results show excellent agreement
  with the high-dimensional scaling of the error predicted by our
  theory.
\end{abstract}

  \end{center}


\section{Introduction}

High-dimensional inference refers to instances of statistical
estimation in which the ambient dimension of the data is comparable to
(or possibly larger than) the sample size. Problems with a
high-dimensional character arise in a variety of applications in
science and engineering, including analysis of gene array data,
medical imaging, remote sensing, and astronomical data analysis.  In
settings where the number of parameters may be large relative to the
sample size, the utility of classical ``fixed $p$'' results is
questionable, and accordingly, a line of on-going statistical research
seeks to obtain results that hold under high-dimensional scaling,
meaning that both the problem size and sample size (as well as other
problem parameters) may tend to infinity simultaneously.  It is
usually impossible to obtain consistent procedures in such settings
without imposing some sort of additional constraints.  Accordingly,
there are now various lines of work on high-dimensional inference
based on imposing different types of structural constraints. A
substantial body of past work has focused on models with sparsity
constraints, including the problem of sparse linear
regression~\cite{Tibshirani96,Chen98,Donoho06,Meinshausen06,BiRiTsy08},
banded or sparse covariance
matrices~\cite{BicLev07,BicLev06,Karoui2007}, sparse inverse
covariance matrices~\cite{YuaLin07,FriedHasTib2007,Rot09,Ravetal08},
sparse eigenstructure~\cite{John01,AmiWai08,PauJoh08}, and sparse
regression matrices~\cite{OboWaiJor08,Lou09,YuaLi06,HuaZha09}.  A
theme common to much of this work is the use of $\ell_1$-penalty as a
surrogate function to enforce the sparsity constraint.

In this paper, we focus on the problem of high-dimensional inference
in the setting of matrix estimation.  As mentioned above, there is
already a substantial body of work on the problem of sparse matrix
recovery.  In contrast, our interest in this paper is the problem of
estimating a matrix $\ThetaStar \in \real^{\lra \times \lrb}$ that is
either \emph{exactly low rank}, meaning that it has at most $\rdim \ll
\min \{\lra, \lrb \}$ non-zero singular values, or more generally is
\emph{near low-rank}, meaning that it can be well-approximated by a
matrix of low rank.  As we discuss at more length in the sequel, such
exact or approximate low-rank conditions are appropriate for many
applications, including multivariate or multi-task forms of
regression, system identification for autoregressive processes,
collaborative filtering, and matrix recovery from random projections.
Analogous to the use of an $\ell_1$-regularizer for enforcing sparsity,
we consider the use of the nuclear norm (also known as the trace norm)
for enforcing a rank constraint in the matrix setting.  By definition,
the nuclear norm is the sum of the singular values of a matrix, and so
encourages sparsity in the vector of singular values, or equivalently
for the matrix to be low-rank.  The problem of low-rank matrix
approximation and the use of nuclear norm regularization have been
studied by various researchers. In her Ph.D. thesis, Fazel~\cite{Fa02}
discusses the use of nuclear norm as a heuristic for restricting the
rank of a matrix, showing that in practice it is often able to yield
low-rank solutions.  Other researchers have provided theoretical
guarantees on the performance of nuclear norm and related methods for
low-rank matrix approximation.  Srebro et al.~\cite{SreRenJaa05}
proposed nuclear norm regularization for the collaborative filtering
problem, and established risk consistency under certain settings.
Recht et al.~\cite{Recht07} provided sufficient conditions for exact
recovery using the nuclear norm heuristic when observing random
projections of a low-rank matrix, a set-up analogous to the compressed
sensing model in sparse linear regression~\cite{Donoho06,CandesTao05}.
Other researchers have studied a version of matrix completion in which
a subset of entries are revealed, and the goal is to obtain perfect
reconstruction either via the nuclear norm heuristic~\cite{CanRec08}
or by other SVD-based methods~\cite{Keshavan08}.  Finally,
Bach~\cite{Bach08} has provided results on the consistency of nuclear
norm minimization for general observation models in noisy settings,
but applicable to the classical ``fixed $p$'' setting.

The goal of this paper is to analyze the nuclear norm relaxation for a
general class of noisy observation models, and obtain non-asymptotic
error bounds on the Frobenius norm that hold under high-dimensional
scaling, and are applicable to both exactly and approximately low-rank
matrices. We begin by presenting a generic observation model, and
illustrating how it can be specialized to the several cases of
interest, including low-rank multivariate regression, estimation of
autoregressive processes, and random projection (compressed sensing)
observations.  In particular, this model is specified in terms of an
operator $\Xop$, which may be deterministic or random depending on the
setting, that maps any matrix $\ThetaStar \in \real^{\lra \times
  \lrb}$ to a vector of $\bigN$ noisy observations.  We then present a
single main theorem (Theorem~\ref{ThmMain}) followed by two
corollaries that cover the cases of exact low-rank constraints
(Corollary~\ref{CorExactLow}) and near low-rank constraints
(Corollary~\ref{CorNearLow}) respectively.  These results demonstrate
that high-dimensional error rates are controlled by two key
quantities.  First, the (random) observation operator $\Xop$ is
required to satisfy a condition known as \emph{restricted strong
  convexity} (RSC), which ensures that the loss function has
sufficient curvature to guarantee consistent recovery of the unknown
matrix $\ThetaStar$.  Second, our theory provides insight into the
\emph{choice of regularization parameter} that weights the nuclear
norm, showing that an appropriate choice is to set it proportional to
the spectral norm of a random matrix defined by the adjoint of
observation operator $\Xop$, and the observation noise in the problem.

This initial set of results, though appealing in terms of their simple
statements and generality, are somewhat abstractly formulated.  Our
next contribution is to show that by specializing our main result
(Theorem~\ref{ThmMain}) to three classes of models, we can obtain some
concrete results based on readily interpretable conditions. In
particular, Corollary~\ref{CorMultivar} deals with the case of
low-rank multivariate regression, relevant for applications in
multitask learning.  We show that the random operator $\Xop$ satisfies
the RSC property for a broad class of observation models, and we use
random matrix theory to provide an appropriate choice of the
regularization parameter.  Our next result,
Corollary~\ref{CorAutoregressive}, deals with the case of estimating
the matrix of parameters specifying a vector autoregressive (VAR)
process~\cite{And71,Lut06}.  Here we also establish that a suitable
RSC property holds with high probability for the random operator
$\Xop$, and also specify a suitable choice of the regularization
parameter.  We note that the technical details here are considerably
more subtle than the case of low-rank multivariate regression, due to
dependencies introduced by the autoregressive sampling scheme.
Accordingly, in addition to terms that involve the size, the matrix
dimensions and rank, our bounds also depend on the mixing rate of the
VAR process.  Finally, we turn to the compressed sensing observation
model for low-rank matrix recovery, as introduced by Recht et
al.~\cite{Recht07}.  In this setting, we again establish that the RSC
property holds with high probability, specify a suitable choice of the
regularization parameter, and thereby obtain a Frobenius error bound
for noisy observations (Corollary~\ref{CorCompressed}).  A technical
result that we prove en route---namely,
Proposition~\ref{LemCompressedRSC}---is of possible independent
interest, since it provides a bound on the constrained norm of a
random Gaussian operator.  In particular, this proposition allows us
to obtain a sharp result (Corollary~\ref{CorRecht}) for the problem of
recovering a low-rank matrix from perfectly observed random
projections, one that removes a logarithmic factor from past
work~\cite{Recht07}.

The remainder of this paper is organized as follows.
Section~\ref{SecBackground} is devoted to background material, and
the set-up of the problem.  We present a generic observation model for
low-rank matrices, and then illustrate how it captures various cases
of interest.  We then define the convex program based on nuclear norm
regularization that we analyze in this paper.  In
Section~\ref{SecResults}, we state our main theoretical results and
discuss their consequences for different model classes.
Section~\ref{SecProofs} is devoted to the proofs of our results; in
each case, we break down the key steps in a series of lemmas, with
more technical details deferred to the appendices.  In
Section~\ref{SecSimulations}, we present the results of various
simulations that illustrate excellent agreement between the
theoretical bounds and empirical behavior.

\noindent \paragraph{Notation:} For the convenience of the reader, we
collect standard pieces of notation here.  For a pair of matrices
$\Theta$ and $\Gamma$ with commensurate dimensions, we let
$\tracer{\Theta}{\Gamma} = \trace(\Theta^T \Gamma)$ denote the trace
inner product on matrix space.  For a matrix $\Theta \in \real^{\lra
  \times \lrb}$, we let $\mymin = \min \{ \lra, \lrb \}$, and denote
its (ordered) singular values by $\sigma_1(\Theta) \geq
\sigma_2(\Theta) \geq \ldots \geq \sigma_\mymin(\Theta) \geq 0$.  We
also use the notation $\sigmax(\Theta) = \sigma_1(\Theta)$ and
$\sigmin(\Theta) = \sigma_{\mymin}(\Theta)$ to refer to the maximal
and minimal singular values respectively.  We use the notation
$\matsnorm{\cdot}{}$ for various types of matrix norms based on these
singular values, including the \emph{nuclear norm} $\nucnorm{\Theta} =
\sum_{j=1}^\mymin \sigma_j(\Theta)$, the \emph{spectral or operator
  norm} $\opnorm{\Theta} = \sigma_1(\Theta)$, and the \emph{Frobenius
  norm} $\matsnorm{\Theta}{F} = \sqrt{\trace(\Theta^T \Theta)} =
\sqrt{\sum_{j=1}^\mymin \sigma_j^2(\Theta)}$.  We refer the reader to
Horn and Johnson~\cite{Horn85,Horn91} for more background on these
matrix norms and their properties.


\section{Background and problem set-up}
\label{SecBackground}

We begin with some background on problems and applications in
which rank constraints arise, before describing a generic observation 
model.  We then introduce the semidefinite program (SDP) based on
nuclear norm regularization that we study in this paper.

\subsection{Models with rank constraints}

Imposing a rank $\rdim$ constraint on a matrix $\ThetaStar \in
\real^{\lra \times \lrb}$ is equivalent to requiring the rows (or
columns) of $\ThetaStar$ lie in some $\rdim$-dimensional subspace of
$\real^{\lrb}$ (or $\real^{\lra}$ respectively).  Such types of rank
constraints (or approximate forms thereof) arise in a variety of
applications, as we discuss here.  In some sense, rank constraints are
a generalization of sparsity constraints; rather than assuming that
the data is sparse in a known basis, a rank constraint implicitly
imposes sparsity but without assuming the basis.

We first consider the problem of multivariate regression, also
referred to as multi-task learning in statistical machine learning.
The goal of \emph{multivariate regression} is to estimate a prediction
function that maps covariates $\Xcov{j} \in \real^\pdim$ to
multi-dimensional output vectors $\Zout{j} \in \real^{\lra}$.  More
specifically, let us consider the linear model, specified by a matrix
$\ThetaStar \in \real^{\lra \times \lrb}$, of the form
\begin{align}
\label{EqnMultivarObs}
\Zout{a} & = \ThetaStar \Xcov{a} + \Wout{a}, \qquad \mbox{for $a = 1,
  \ldots, \numobs$},
\end{align}
where $\{ \Wout{a} \}_{a=1}^\numobs$ is an i.i.d. sequence of
$\lra$-dimensional zero-mean noise vectors.  Given a collection of
observations $\{\Xcov{a}, \Zout{a}\}_{a=1}^\numobs$ of
covariate-output pairs, our goal is to estimate the unknown matrix
$\ThetaStar$.  This type of model has been used in many applications,
including analysis of fMRI image data~\cite{HarPenFri03}, analysis of
EEG data decoding~\cite{An98}, neural response
modeling~\cite{BroKasMit04} and analysis of financial data.  This
model and closely related ones also arise in the problem of
collaborative filtering~\cite{SreRenJaa04}, in which the goal is to
predict users' preferences for items (such as movies or music) based
on their and other users' ratings of related items.  The
papers~\cite{Abeplus06,ArgEvgPon06} discuss additional instances of
low-rank decompositions.  In all of these settings, the low-rank
condition translates into the existence of a smaller set of
``features'' that are actually controlling the prediction.

As a second (not unrelated) example, we now consider the problem of
system identification in vector autoregressive processes (see the
book~\cite{Lut06} for detailed background).  A \emph{vector
autoregressive} (VAR) process in $\pdim$-dimensions is a a stochastic
process $\{\Zsam{t}\}_{t=1}^\infty$ specified by an initialization
$\Zsam{1} \in \real^\pdim$, followed by the recursion
\begin{align}
\label{EqnAuto}
  \Zsam{t+1} & = \ThetaStar \Zsam{t} + \Noisesam{t}, \qquad \mbox{for
    $t = 1, 2, 3, \ldots$}.
\end{align}
In this recursion, the sequence $\{\Noisesam{t}\}_{t=1}^\infty$
consists of i.i.d. samples of innovations noise.  We assume that each
vector $\Noisesam{t} \in \real^\pdim$ is zero-mean with covariance
$\nvar^2 I$, so that the process $\{\Zsam{t}\}_{t=1}^\infty$ is
zero-mean, and has a covariance matrix $\CovMat$ given by the solution
of the discrete-time Ricatti equation
\begin{align}
\label{EqnRicatti}
\CovMat & = \ThetaStar \CovMat (\ThetaStar)^T + \nvar^2 I.
\end{align}
The goal of system identification in a VAR process is to estimate the
unknown matrix \mbox{$\ThetaStar \in \real^{\pdim \times \pdim}$} on
the basis of a sequence of samples
\mbox{$\{\Zsam{t}\}_{t=1}^\numobs$.}  In many application domains, it
is natural to expect that the system is controlled primarily by a
low-dimensional subset of variables.  For instance, models of
financial data might have an ambient dimension $\pdim$ of thousands
(including stocks, bonds, and other financial instruments), but the
behavior of the market might be governed by a much smaller set of
macro-variables (combinations of these financial instruments).
Similar statements apply to other types of time series data, including
neural data~\cite{BroKasMit04,FisBla05}, subspace tracking models in
signal processing, and motion models models in computer vision.

A third example that we consider in this paper is a \emph{compressed
sensing} observation model, in which one observes random projections
of the unknown matrix $\ThetaStar$.  This observation model has been
studied extensively in the context of estimating sparse
vectors~\cite{Donoho06,CandesTao05}, and Recht et al.~\cite{Recht07}
suggested and studied its extension to low-rank matrices.  In their
set-up, one observes trace inner products of the form
$\tracer{X_i}{\ThetaStar} = \trace(X_i^T \ThetaStar)$, where $X_i \in
\real^{\lra \times \lrb}$ is a random matrix (for instance, filled
with standard normal $N(0,1)$ entries).  Like compressed sensing for
sparse vectors, applications of this model include computationally
efficient updating in large databases (where the matrix $\ThetaStar$
measures the difference between the data base at two different time
instants), and matrix denoising.


\subsection{A generic observation model}

We now introduce a generic observation model that will allow us to
deal with these different observation models in an unified manner.
For pairs of matrices $A, B \in \real^{\lra \times \lrb}$, recall the
Frobenius or trace inner product $\tracer{A}{B} \defn \trace(B A^T)$.
We then consider a linear observation model of the form
\begin{align}
\label{EqnSmallObs}
\mysmally{i} & = \tracer{\mysmallX{i}}{\ThetaStar} + \mysmallw{i},
\qquad \mbox{for $i=1, 2, \ldots, \bigN$,}
\end{align}
which is specified by the sequence of observation matrices
$\{\mysmallX{i}\}_{i=1}^\bigN$ and observation noise
$\{\mysmallw{i}\}_{i=1}^\bigN$.
This observation model can be written in a more compact manner using
operator-theoretic notation.  In particular, let us define the
observation vector 
\begin{align*}
\mybigY  & = 
\begin{bmatrix} \mysmally{1} & \ldots & \mysmally{\numobs}
\end{bmatrix}^T \in \real^\bigN,
\end{align*}
with a similar definition for $\mybigW \in \real^\bigN$ in terms of
$\{\mysmallw{i}\}_{i=1}^\bigN$.  We then use the observation matrices
$\{\mysmallX{i}\}_{i=1}^\bigN$ to define an operator \mbox{$\Xop:
  \real^{\lra \times \lrb} \rightarrow \real^\bigN$} via
$\big[\Xop(\Theta)\big]_i = \tracer{\mysmallX{i}}{\Theta}$.  With this
notation, the observation model~\eqref{EqnSmallObs} can be re-written as
\begin{align}
\label{EqnCompactObs}
\mybigY & = \Xop(\ThetaStar) + \mybigW.
\end{align}

Let us illustrate the form of the observation
model~\eqref{EqnCompactObs} for some of the applications that we
considered earlier.

\begin{exas}[Multivariate regression]
\label{ExaMultivar}
Recall the observation model~\eqref{EqnMultivarObs} for multivariate
regression.  In this case, we make $\numobs$ observations of vector
pairs \mbox{$(\Zout{a}, \Xcov{a}) \in \real^\lra \times
\real^{\lrb}$.}  Accounting for the $\lra$-dimensional nature of the
output, after the model is scalarized, we receive a total of $\bigN =
\lra \numobs$ observations.  Let us introduce the quantity $b = 1,
\ldots, \lra$ to index the different elements of the output, so that
we can write
\begin{align}
\Zout{ab} & = \tracer{\Xcov{a} e_b^T}{\ThetaStar} + \Wout{ab}, \qquad
\mbox{for $b = 1, 2, \ldots, \lra$.}
\end{align}
By re-indexing this collection of $\bigN = \numobs \lra$ observations
via the mapping $(a,b) \mapsto i = a + (b-1) \, \lra$, we recognize
multivariate regression as an instance of the observation
model~\eqref{EqnSmallObs} with observation matrix $\mysmallX{i} =
\Xcov{a} e_b^T$ and scalar observation $y_i = V_{ab}$.
\end{exas}


\begin{exas}[Vector autoregressive processes]
\label{ExaAutoregressive}
Recall that a vector autoregressive (VAR) process is defined by the
recursion~\eqref{EqnAuto}, and suppose that observe an
$\numobs$-sequence $\{\Zsam{t}\}_{t=1}^\numobs$ produced by this
recursion.  Since each $\Zsam{t} =
\begin{bmatrix} \Zsam{t1} & \ldots & \Zsam{t\pdim} \end{bmatrix}^T$ is
$\pdim$-variate, the scalarized sample size is $\bigN = \numobs
\pdim$.  Letting $b = 1, 2, \ldots, \pdim$ index the dimension, we
have
\begin{align}
\Zsam{(t+1) \, b} & = \tracer{\Zsam{t}e_b^T}{\ThetaStar} + \Noisesam{tb}.
\end{align}
In this case, we re-index the collection of $\bigN = \numobs \pdim$
observations via the mapping $(t,b) \mapsto i = t + (b-1) \, \pdim$.
After doing so,  we see that the autoregressive problem can be written
in the form~\eqref{EqnSmallObs} with $y_i = \Zsam{(t+1) \, b}$ and
observation matrix $\mysmallX{i} = \Zsam{t} e_b^T$.
\end{exas}


\begin{exas}[Compressed sensing]
\label{ExaCompressed}
As mentioned earlier, this is a natural extension of the compressed
sensing observation model for sparse vectors to the case of low-rank
matrices~\cite{Recht07}.  In particular, suppose that each observation
matrix $\mysmallX{i} \in \real^{\lra \times \lrb}$ has i.i.d. standard
normal $N(0,1)$ entries, so that we make observations of the form
\begin{align}
y_i & = \tracer{\mysmallX{i}}{\ThetaStar} + \mysmallw{i}, \quad
\mbox{for $i = 1, 2, \ldots, \numobs$.}
\end{align}
By construction, these observations are an instance of the
model~\eqref{EqnSmallObs}.  In this case, the more compact
form~\eqref{EqnCompactObs} involves a random Gaussian operator mapping
$\real^{\lra \times \lrb}$ to $\real^\numobs$, and we study some of its
properties in the sequel.
\end{exas}

\subsection{Regression with nuclear norm regularization}

We now consider an estimator that is naturally suited to the problems
described in the previous section.  Recall that the \emph{nuclear or
trace norm} of a matrix $\Theta \in \real^{\lra \times \lrb}$ is given
by $\nucnorm{\Theta} = \sum_{j=1}^{\min\{\lra, \lrb\}}
\sigma_j(\Theta)$, corresponding to the sum of its singular values.
Given a collection of observations $(\mysmally{i}, \mysmallX{i}) \in
\real \times \real^{\lra \times \lrb}$, for $i = 1, \ldots, \bigN$
from the observation model~\eqref{EqnSmallObs}, we consider estimating
the unknown $\ThetaStar$ by solving the following optimization problem
\begin{align}
\label{EqnSDP}
\ThetaHat & \in \arg \min_{\Theta \in \real^{\lra \times \lrb}} \Big
\{ \frac{1}{2 \bigN} \|\mybigY - \Xop(\Theta) \|_2^2 + \bigreg
\matsnorm{\Theta}{1} \Big \},
\end{align}
where $\bigreg > 0$ is a regularization parameter.  Note that the
optimization problem~\eqref{EqnSDP} can be viewed as the analog of the
Lasso estimator~\cite{Tibshirani96}, tailored to low-rank matrices as
opposed to sparse vectors.  An important property of the optimization
problem~\eqref{EqnSDP} is that it can be solved in time polynomial in
the sample size $\bigN$ and the matrix dimensions $\lra$ and $\lrb$.
Indeed, the optimization problem~\eqref{EqnSDP} is an instance of a
\emph{semidefinite program}~\cite{Vanden96}, a class of convex
optimization problems that can be solved efficiently by various
polynomial-time algorithms~\cite{Boyd02}.  For instance, interior
point methods are a classical method for solving semidefinite
programs; moreover, as we discuss in Section~\ref{SecSimulations},
there are a variety of other methods for solving the semidefinite
program (SDP) defining our $M$-estimator

Like in any typical $M$-estimator for statistical inference, the
regularization parameter $\bigreg$ is specified by the statistician.  As
part of the theoretical results in the next section, we provide
suitable choices of this parameter in order for the estimate
$\ThetaHat$ to behave well, in the sense of being close in Frobenius
norm to the unknown matrix $\ThetaStar$.


\section{Main results and some consequences}
\label{SecResults}

In this section, we state our main results and discuss some of their
consequences.  Section~\ref{SecResultsGeneric} is devoted to results
that apply to generic instances of low-rank problems, whereas
Section~\ref{SecResultsSpecific} is devoted to the consequences of
these results for more specific problem classes, including low-rank
multivariate regression, estimation of vector autoregressive
processes, and recovery of low-rank matrices from random projections.

\subsection{Results for general model classes}
\label{SecResultsGeneric}

We begin by introducing the key technical condition that allows us to
control the error $\ThetaHat - \ThetaStar$ between an SDP solution
$\ThetaHat$ and the unknown matrix $\ThetaStar$.  We refer to it as
the \emph{restricted strong convexity} condition~\cite{NegRavWaiYu09},
since it amounts to guaranteeing that the quadratic loss function in
the SDP~\eqref{EqnSDP} is strictly convex over a restricted set of
directions.  Letting $\KeySet \subseteq \real^{\lra \times \lrb}$
denote the restricted set of directions, we say that the operator
$\Xop$ satisfies restricted strong convexity (RSC) over the set
$\KeySet$ if there exists some $\kappa(\Xop) > 0$ such that
\begin{align}
\label{EqnDefnRSC}
\frac{1}{2 \bigN} \|\Xop(\Delta)\|_2^2 & \geq \kappa(\Xop) \;
\matsnorm{\Delta}{F}^2 \qquad \mbox{for all $\Delta \in \KeySet$.}
\end{align}
We note that analogous conditions have been used to establish error
bounds in the context of sparse linear
regression~\cite{BiRiTsy08,CohDahDev08}, in which case the set
$\KeySet$ corresponded to certain subsets of sparse vectors.

Of course, the definition~\eqref{EqnDefnRSC} hinges on the choice of
the restricted set $\KeySet$.  In order to specify some appropriate
sets for the case of low-rank matrices, we require some additional
notation.  For any matrix $\Theta \in \real^{\lra \times \lrb}$, we
let $\row(\Theta) \subseteq \real^{\lrb}$ and $\col(\Theta) \subseteq
\real^{\lra}$ denote its row space and column space, respectively.
For a given positive integer $\rdim \leq \min \{ \lra, \lrb \}$, any
$\rdim$-dimensional subspace of $\real^\lra$ can be represented by
some orthogonal matrix $U \in \real^{\lra \times \rdim}$ (i.e., that
satisfies $U^T U = I_{\rdim \times \rdim}$.  In a similar fashion, any
$\rdim$-dimensional subspace of $\real^{\lrb}$ can be represented by
an orthogonal matrix $V \in \real^{\lrb \times \rdim}$.  For any fixed
pair of such matrices $(\mata, \matb)$, we may define the following
two subspaces of $\real^{\lra \times \lrb}$:
\begin{subequations}
\label{EqnDefnSetOneTwo}
\begin{align}
\SetOne(\mata, \matb) & \defn \big \{ \Theta \in \real^{\lra \times
  \lrb} \, \mid \, \mbox{$\row(\Theta) = \matb$ and $\col(\Theta) =
  \mata$} \big \}, \qquad \mbox{and} \\
\SetTwo(\mata, \matb) & \defn \big \{ \Theta \in \real^{\lra \times
  \lrb} \, \mid \, \mbox{$\row(\Theta) \perp \matb$ and
  $\col(\Theta) \perp  \mata$.} \big \}.
\end{align}
\end{subequations}
Finally, we let $\Pi_{\SetOne(\mata, \matb)}$ and $\Pi_{\SetTwo(\mata,
\matb)}$ denote the (respective) projection operators onto these
subspaces.  When the subspaces $(\mata, \matb)$ are clear from
context, we use the shorthand notation $\PlErrorB =
\Pi_{\SetTwo(\mata, \matb)}(\Error)$ and $\PlErrorA = \Delta -
\PlErrorB$.  With this notation, we can define the restricted sets
$\KeySet$ of interest.  Using the singular value decomposition, we can
write \mbox{$\ThetaStar = \mataspec D \matbspec^T$,} where $\mataspec
\in \real^{\lra \times \lra}$ and $\matbspec \in \real^{\lrb \times
\lrb}$ are both orthogonal matrices, and $D \in \real^{\lra \times
\lrb}$ contains the singular values of $\ThetaStar$.  For any positive
integer $\rdim \leq \min \{\lra, \lrb\}$, we let $(\matardim,
\matbrdim)$ denote the subspace pair defined by the top $\rdim$ left
and right singular vectors of $\ThetaStar$.  For a given integer
$\rdim$ and tolerance $\delta > 0$, we then define a subset of
matrices as follows:
\begin{align}
\label{eq:restrict}
\KeySet(\rdim; \delta) & \mydefn \biggr \{ \Delta \in \real^{\lra
\times \lrb} \, \mid \, \matsnorm{\Delta}{F} \geq \delta, \;
\nucnorm{\PlErrorB} \; \leq \; 3 \nucnorm{\PlErrorA} \, + \, 4
\matsnorm{\Pi_{\SetTwo(\matardim, \matbrdim)}(\ThetaStar)}{1} \biggr
\}.
\end{align}

The next ingredient is the choice of the regularization parameter
$\bigreg$ used in solving the SDP~\eqref{EqnSDP}.  Our theory
specifies a choice for this quantity in terms of the adjoint of the
operator $\Xop$---namely, the operator $\Xopadj:\real^\bigN
\rightarrow \real^{\lra \times \lrb}$ defined by
\begin{align}
\label{EqnDefnXopadj}
\Xopadj(\mybigW) & \mydefn \sum_{i=1}^{\bigN} \mysmallw{i} X_i.
\end{align}
With this notation, we come to the first result of our paper.  It is a
deterministic result, which specifies two conditions---namely, an RSC
condition and a choice of the regularizer---that suffice to guarantee
for any solution of the SDP~\eqref{EqnSDP} fall within a certain
radius.
\btheos
\label{ThmMain}
Suppose that the operator $\Xop$ satisfies restricted strong convexity
with parameter $\kappa(\Xop) > 0$ over the set $\KeySet(\rdim;
\delta)$, and that the regularization parameter $\bigreg$ is chosen
such that \mbox{$\bigreg \geq 2
\opnorm{\Xopadj(\mybigW)}/\bigN$.} Then any solution $\ThetaHat$
to the semidefinite program~\eqref{EqnSDP} satisfies
\begin{align}
\label{EqnGenError}
\matsnorm{\ThetaHat - \ThetaStar}{F} & \leq \max \biggr \{\delta, \;
\frac{32 \bigreg \: \sqrt{\rdim}}{\kappa(\Xop)}, \; \biggr[\frac{16 \:
\bigreg \: \matsnorm{\Pi_{
\SetTwo(\matardim,
\matbrdim)}(\ThetaStar)}{1}}{\kappa(\Xop)} \biggr]^{1/2} \biggr \}.
\end{align}
\etheos

Apart from the tolerance parameter $\delta$, the two main terms in the
bound~\eqref{EqnGenError} have a natural interpretation. The first
term (involving $\sqrt{\rdim}$) corresponds to \emph{estimation
error}, capturing the difficulty of estimating a rank $\rdim$ matrix.
The second is an \emph{approximation error}, in which the projection
onto the set $\SetTwo(\matardim, \matbrdim)$ describes the gap between
the true matrix $\ThetaStar$ and the rank $\rdim$ approximation. \\

Let us begin by illustrating the consequences of Theorem~\ref{ThmMain}
when the true matrix $\ThetaStar$ has exactly rank $\rdim$, in which case
there is a very natural choice of the subspaces represented by $\mata$
and $\matb$.  In particular, we form $\mata$ from the $\rdim$ non-zero
left singular vectors of $\ThetaStar$, and $\matb$ from its $\rdim$
non-zero right singular vectors.  Note that this choice of $(\mata,
\matb)$ ensures that $\Pi_{\SetTwo(\mata, \matb)}(\ThetaStar) = 0$.
For technical reasons to be clarified, it suffices to set $\delta = 0$
in the case of exact rank constraints, and we thus obtain the
following result:
\bcors[Exact low-rank recovery]
\label{CorExactLow}
Suppose that $\ThetaStar$ has rank $\rdim$, and $\Xop$
satisfies RSC with respect to $\KeySet(\rdim; 0)$.  Then as long as
\mbox{$\bigreg \geq 2 \opnorm{\Xopadj(\mybigW)}/\bigN$,} any
optimal solution $\ThetaHat$ to the SDP~\eqref{EqnSDP} satisfies the
bound
\begin{align}
\label{EqnLowRankError}
\matsnorm{\ThetaHat - \ThetaStar}{F} & \leq \frac{32 \sqrt{\rdim} \;
\bigreg}{\kappa(\Xop)}.
\end{align}
\ecors

\noindent Like Theorem~\ref{ThmMain}, Corollary~\ref{CorExactLow} is a
deterministic statement on the SDP error.  It takes a much simpler
form since when $\ThetaStar$ is exactly low rank, then neither
tolerance parameter $\delta$ nor the approximation term are required. \\

As a more delicate example, suppose instead that $\ThetaStar$ is
\emph{nearly low-rank}, an assumption that we can formalize by
requiring that its singular value sequence $\{\sigma_i(\ThetaStar)
\}_{i=1}^{\min\{ \lra, \lrb \}}$ decays quickly enough.  In particular,
for a parameter $\qpar \in [0,1]$ and a positive radius $\radq$, we
define the set
\begin{align}
\weakball & \defn \big \{ \Theta \in \real^{\lra \times \lrb} \, \mid
\, \sum_{i=1}^{\min \{ \lra, \lrb \}} |\sigma_i(\Theta)|^\qpar \leq
\radq \big \}.
\end{align}
Note that when $\qpar = 0$, the set $\ellzball$ corresponds to
the set of matrices with rank at most $\radz$.

\bcors[Near low-rank recovery]
\label{CorNearLow}
Suppose that $\ThetaStar \in \weakball$, the regularization parameter
is lower bounded as \mbox{$\bigreg \geq 2
\opnorm{\Xopadj(\mybigW)}/\bigN$,} and the operator $\Xop$
satisfies RSC with parameter $\kappa(\Xop) \in (0,1]$ over the set
$\KeySet(\radq \bigreg^{-\qpar}; \delta)$.  Then any solution
$\ThetaHat$ to the SDP~\eqref{EqnSDP} satisfies
\begin{align}
\label{EqnNearLowError}
\matsnorm{\ThetaHat - \ThetaStar}{F} & \leq \max \Big \{ \delta, \; 32
\; \sqrt{\radq} \biggr (\frac{\bigreg}{\kappa(\Xop)}
\biggr)^{1-\qpar/2} \Big \}.
\end{align}
\ecors
\noindent Note that the error bound~\eqref{EqnNearLowError} reduces to
the exact low rank case~\eqref{EqnLowRankError} when $\qpar = 0$, and
$\delta = 0$.  The quantity $\bigreg^{-\qpar} \radq$ acts as the
``effective rank'' in this setting; as clarified by our proof in
Section~\ref{SecProofCorNearLow}. This particular choice is designed
to provide an optimal trade-off between the approximation and
estimation error terms in Theorem~\ref{ThmMain}.  Since $\bigreg$ is
chosen to decay to zero as the sample size $\bigN$ increases, this
effective rank will increase, reflecting the fact that as we obtain
more samples, we can afford to estimate more of the smaller singular
values of the matrix $\ThetaStar$.

\subsection{Results for specific model classes}
\label{SecResultsSpecific}

As stated, Corollaries~\ref{CorExactLow} and~\ref{CorNearLow} are
fairly abstract in nature.  More importantly, it is not immediately
clear how the key underlying assumption---namely, the RSC
condition---can be verified, since it is specified via subspaces that
depend on $\ThetaStar$, which is itself the unknown quantity that we
are trying to estimate.  Nonetheless, we now show how, when
specialized to more concrete models, these results yield concrete and
readily interpretable results.  As will be clear in the proofs of
these results, each corollary requires overcoming two main technical
obstacles: establishing that the appropriate form of the RSC property
holds in a uniform sense (so that a priori knowledge of $\ThetaStar$
is not required), and specifying an appropriate choice of the
regularization parameter $\bigreg$.  Each of these two steps is
non-trivial, requiring some random matrix theory, but the end results
are simply stated upper bounds that hold with high probability.

We begin with the case of rank-constrained multivariate regression.
As discussed earlier in Example~\ref{ExaMultivar}, recall that we
observe pairs $(\Zout{i}, \Xcov{i}) \in \real^{\lra} \times
\real^{\lrb}$ linked by the linear model $\Zout{i} = \ThetaStar
\Xcov{i} + \Wout{i}$, where $\Wout{i} \sim N(0, \nvar^2 I_{\lra \times
\lra})$ is observation noise.  Here we treat the case of \emph{random
design regression}, meaning that the covariates $\Xcov{i}$ are modeled
as random.  In particular, in the following result, we assume that
$\Xcov{i} \sim N(0, \Sigma)$, i.i.d.  for some $\pdim$-dimensional
covariance matrix $\Sigma \succ 0$.  Recalling that $\sigmax(\Sigma)$
and $\sigmin(\Sigma)$ denote the maximum and minimum eigenvalues
respectively, we have:
\bcors[Low-rank multivariate regression]
\label{CorMultivar}
Consider the random design multivariate regression model where
$\ThetaStar \in \weakball$.  There are universal constants $\{c_i,
i=1,2,3 \}$ such that if we solve the SDP~\eqref{EqnSDP} with
regularization parameter $\bigreg = 10 \nvar \sqrt{\sigmax(\Sigma)} \;
\sqrt{\frac{(\lra + \lrb)}{\numobs}}$, we have
\begin{align}
\label{EqnMultivarError}
\matsnorm{\ThetaHat - \ThetaStar}{F}^2 & \leq c_1 \;
\biggr(\frac{\nvar^2 \sigmax(\Sigma)}{\sigmin^2(\Sigma)}
\biggr)^{1-\qpar/2} \; \radq \; \biggr(\frac{\lra + \lrb}{\numobs}
\biggr)^{1-\qpar/2}
\end{align}
with probability greater than $1 - c_2 \exp(-c_3 (\lra + \lrb))$.
\ecors

\noindent {\bf{Remarks:}} Corollary~\ref{CorMultivar} takes a
particularly simple form when $\Sigma = I_{\pdim \times \pdim}$: then
there exists a constant $c'_1$ such that $\matsnorm{\ThetaHat -
\ThetaStar}{F}^2 \leq c'_1 \nvar^2 \: \radq \; \big(\frac{\lra +
\lrb}{\numobs}\big)^{1-\qpar/2}$.  When $\ThetaStar$ is exactly low
rank---that is, $\qpar = 0$, and $\ThetaStar$ has rank $\rdim =
R_0$---this simplifies even further to
\begin{align*}
\matsnorm{\ThetaHat - \ThetaStar}{F}^2 \leq c'_1 \frac{ \nvar^2 \;
\rdim \, (\lra + \lrb)}{\numobs}.
\end{align*}
The scaling in this error bound is easily interpretable: naturally,
the squared error is proportional to the noise variance $\nvar^2$, and
the quantity $\rdim (\lra + \lrb)$ counts the number of degrees of
freedom of a $\lra \times \lrb$ matrix with rank $\rdim$.  Note that
if we did not impose any constraints on $\ThetaStar$, then since a
$\lra \times \lrb$ matrix has a total of $\lra \lrb$ free parameters,
we would expect at best to obtain rates of the order
$\matsnorm{\ThetaHat - \ThetaStar}{F}^2 = \Omega(\frac{\nvar^2 \, \lra
\, \lrb}{\numobs})$.  Note that when $\ThetaStar$ is low rank---in
particular, when $\rdim \ll \min \{\lra,\lrb \}$---then the nuclear
norm estimator achieves substantially faster rates.  Finally, we note
that as stated, the result requires that $\min \{\lra, \lrb \}$ tend
to infinity in order for the claim to hold with high probability.
Although such high-dimensional scaling is the primary focus of this
paper, we note that for application to the classical setting of fixed
$(\lra, \lrb)$, the same statement (with different constants) holds
with $\lra + \lrb$ replaced by $\log \numobs$. \\

Next we turn to the case of estimating the system matrix $\ThetaStar$
of an autoregressive (AR) model, as discussed in
Example~\ref{ExaAutoregressive}.

\bcors[Autoregressive models] 
\label{CorAutoregressive}
Suppose that we are given $\numobs$ samples
$\{\Zsam{t}\}_{t=1}^\numobs$ from a $\pdim$-dimensional autoregressive
process~\eqref{EqnAuto} that is stationary, based on a system matrix
that is stable ($\opnorm{\ThetaStar} \leq \maxop < 1$), and
approximately low-rank ($\ThetaStar \in \weakball$).  Then there are
universal constants $\{c_i, i=1,2,3 \}$ such that if we solve the
SDP~\eqref{EqnSDP} with regularization parameter $\bigreg = \frac{80
  \, \opnorm{\Sigma}}{1-\maxop} \sqrt{\frac{\pdim}{\numobs}}$, then
any solution $\ThetaHat$ satisfies
\begin{align}
\label{EqnAutoregressiveError}
  \frobnorm{\widehat{\paramvar} - \thetastar}^2 & \leq c_1 \; \biggr[
  \frac{\sigmax(\Sigma)}{\sigmin(\Sigma) \, (1 - \maxop)}
  \biggr]^{2-\qpar} \, \radq \; \Big(\frac{\pdim}{\numobs}
  \Big)^{1-\qpar/2}
\end{align}
with probability greater than $1-c_2 \exp(-c_3 \pdim)$.
\ecors

\noindent {\bf{Remarks:}} Like Corollary~\ref{CorMultivar}, the result
as stated requires that $\pdim$ tend to infinity, but the same bounds
hold with $\pdim$ replaced by $\log \numobs$, yielding results
suitable for classical (fixed dimension) scaling.  Second, the factor
$(\pdim/\numobs)^{1-\qpar/2}$, like the analogous term\footnote{The
term in Corollary~\ref{CorMultivar} has a factor $\lra + \lrb$, since
the matrix in that case could be non-square in general.}  in
Corollary~\ref{CorMultivar}, shows that faster rates are obtained if
$\ThetaStar$ can be well-approximated by a low rank matrix, namely for
choices of the parameter $\qpar \in [0,1]$ that are closer to zero.
Indeed, in the limit $\qpar = 0$, we again reduce to the case of an
exact rank constraint $\rdim = R_0$, and the corresponding squared
error scales as $\rdim \pdim/\numobs$.  In contrast to the case of
multivariate regression, the error
bound~\eqref{EqnAutoregressiveError} also depends on the upper bound
$\opnorm{\ThetaStar} = \maxop < 1$ on the operator norm of the system
matrix $\ThetaStar$.  Such dependence is to be expected since the
quantity $\maxop$ controls the (in)stability and mixing rate of the
autoregressive process. As clarified in the proof, the dependence of
the sampling in the AR model also presents some technical challenges
not present in the setting of multivariate regression. \\

Finally, we turn to analysis of the compressed sensing model for
matrix recovery, which was discussed in Example~\ref{ExaCompressed}.
The following result applies to the setting in which the observation
matrices $\{\mybigX_i\}_{i=1}^\bigN$ are drawn i.i.d., with standard
$N(0,1)$ elements.  We assume that the observation noise vector
$\mybigW \in \real^\bigN$ satisfies the bound $\|\mybigW\|_2 \leq 2
\nvar \, \sqrt{\bigN}$ for some constant $\nvar$, an assumption that
holds for any bounded noise, and also holds with high probability for
any random noise vector with sub-Gaussian entries with parameter
$\nvar$ (one example being Gaussian noise $N(0, \nvar^2)$).

\bcors[Compressed sensing recovery]
\label{CorCompressed}
Suppose that $\ThetaStar \in \weakball$, and that the sample size is
lower bounded as $\bigN > 200 \, \radq^{1-\qpar/2} \; \lra \, \lrb$.
Then any solution $\ThetaHat$ to the SDP~\eqref{EqnSDP} satisfies the
bound
\begin{align}
  \matsnorm{\ThetaHat - \ThetaStar}{F}^2 & \leq 256 \; \nvar
^{2-\qpar} \; \radq \; \Big [ \sqrt{\frac{\lra}{\bigN}} +
\sqrt{\frac{\lrb}{\bigN}} \Big ]^{2 - \qpar}
\end{align}
with probability greater than $1 - c_1 \exp(-c_2 (\lra + \lrb))$.
\ecors

The central challenge in proving this result is in proving an
appropriate form of the RSC property.  The following result on the
random operator $\Xop$ may be of independent interest here:

\bprops
\label{LemCompressedRSC}
Under the stated conditions, the random operator $\Xop$ satisfies
\begin{align}
\label{EqnCompressedRSC}
\frac{\| \Xop(\Theta)\|_2}{\sqrt{\bigN}} & \geq \conone
\matsnorm{\Theta}{F} - \contwo \BigMess
\nucnorm{\Theta} \qquad \mbox{for all $\Theta \in \real^{\lra \times
    \lrb}$}
\end{align}
with probability at least $1 - 2\exp(-\bigN/32)$.
\eprops
\noindent The proof of this result, provided in
Appendix~\ref{AppLemCompressedRSC}, makes use of the Gordon-Slepian
inequalities for Gaussian processes, and concentration of measure.  As
we show in Section~\ref{SecProofCorCompressed}, it implies the form of
the RSC property needed to establish Corollary~\ref{CorCompressed}.

Proposition~\ref{LemCompressedRSC} also implies an interesting
property of the null space of the operator $\Xop$; one that can be
used to establish a corollary about recovery of low-rank matrices when
the observations are noiseless.  In particular, suppose that we are
given the noiseless observations $y_i = \tracer{X_i}{\ThetaStar}$ for
$i = 1, \ldots, \bigN$, and that we try to recover the unknown matrix
$\ThetaStar$ by solving the SDP
\begin{equation}
\label{EqnBasisPursuitSDP}
\min_{\Theta \in \real^{\lra \ltimes \lrb}} \nucnorm{\Theta} \quad
\mbox{such that $\tracer{X_i}{\Theta} = y_i$ for all $i = 1, \ldots,
\bigN$,}
\end{equation}
a recovery procedure that was studied by Recht et al.~\cite{Recht07}.
Proposition~\ref{LemCompressedRSC} allows us to obtain a sharp result
on recovery using this method: \\

\bcors
\label{CorRecht}
Suppose that $\ThetaStar$ has rank $\rdim$, and that we are given
$\bigN > 40^2 \rdim (\lra + \lrb)$ noiseless samples.  Then with
probability at least $1 - 2\exp(-\bigN/32)$, the
SDP~\eqref{EqnBasisPursuitSDP} recovers the matrix $\ThetaStar$
exactly.
\ecors
\noindent This result removes some extra logarithmic factors that were
included in the earlier work~\cite{Recht07}, and provides the
appropriate analog to compressed sensing results for sparse
vectors~\cite{Donoho06,CandesTao05}.  Note that (like in most of our
results) we have made little effort to obtain good constants in this
result: the important property is that the sample size $\bigN$ scales
linearly in both $\rdim$ and $\lra + \lrb$.


\section{Proofs}
\label{SecProofs}

We now turn to the proofs of Theorem~\ref{ThmMain}, and
Corollaries~\ref{CorExactLow} through~\ref{CorRecht}.  In each
case, we provide the primary steps in the main text, with more
technical details stated as lemmas and proved in the Appendix.

\subsection{Proof of Theorem~\ref{ThmMain}}

By the optimality of $\ThetaHat$ for the SDP~\eqref{EqnSDP}, we have
\begin{align*}
\frac{1}{2 \bigN} \|\mybigY - \Xop(\ThetaHat) \|_2^2 + \bigreg
\matsnorm{\ThetaHat}{1} & \leq \frac{1}{2 \bigN} \|\mybigY -
\Xop(\ThetaStar) \|_2^2 + \bigreg \matsnorm{\ThetaStar}{1}.
\end{align*}
Defining the error matrix $\PlError = \ThetaStar - \ThetaHat$ and
performing some algebra yields the inequality
\begin{align}
\label{EqnBasic}
\frac{1}{2 \bigN} \|\Xop(\Delta)\|_2^2 & \leq \frac{1}{\bigN}
\inprod{\mybigW}{\Xop(\PlError)} + \bigreg \big \{
\matsnorm{\ThetaHat + \PlError }{1} - \matsnorm{\ThetaHat}{1} \big \}.
\end{align}
By definition of the adjoint and H\"{o}lder's inequality, we have
\begin{align}
\label{EqnCafeOne}
\frac{1}{\bigN} \big|\inprod{\mybigW}{\Xop(\PlError)} \big| & =
\frac{1}{\bigN} \big|\inprod{\Xopadj(\mybigW)}{\PlError} \big| \;  \leq
\frac{1}{\bigN} \matsnorm{\Xopadj(\mybigW)}{\myop} \,
\matsnorm{\PlError}{1}.
\end{align}
By the triangle inequality, we have $\matsnorm{\ThetaHat + \PlError
}{1} - \matsnorm{\ThetaHat}{1} \leq \matsnorm{\PlError}{1}$.
Substituting this inequality and the bound~\eqref{EqnCafeOne} into the
inequality~\eqref{EqnBasic} yields
\begin{align*}
\frac{1}{2 \bigN} \|\Xop(\Delta)\|_2^2 & \leq \frac{1}{\bigN}
\matsnorm{\Xopadj(\mybigW)}{\myop} \matsnorm{\PlError}{1} +\bigreg
\matsnorm{\PlError}{1} \; \leq \; 2 \bigreg \matsnorm{\PlError}{1},
\end{align*}
where the second inequality makes use of our choice $\bigreg \geq
\frac{2}{\bigN} \matsnorm{\Xopadj(\mybigW)}{\myop}$.

It remains to lower bound the term on the left-hand side, while upper
bounding the quantity $\matsnorm{\PlError}{1}$ on the right-hand side.
The following technical result allows us to do so.  Recall our earlier
definition~\eqref{EqnDefnSetOneTwo} of the sets $\SetOne$ and
$\SetTwo$ associated with a given subspace pair.

\blems
\label{LemRestricted}
Let $(\mataspec, \matbspec)$ represent a pair of $\rdim$-dimensional
subspaces of left and right singular vectors of $\ThetaStar$.  Then
there exists a matrix decomposition $\PlError = \PlErrorA + \PlErrorB$
of the error $\PlError$ such that
\begin{enumerate}
\item[(a)] The matrix $\PlErrorA$ satisfies the constraint
$\rank(\PlErrorA) \leq 2 \rdim$, and
\item[(b)] If \mbox{$\bigreg \geq 2
\matsnorm{\Xopadj(\mybigW)}{2}/\bigN$,} then the nuclear norm of
$\PlErrorB$ is bounded as
\begin{equation}
\label{EqnRestricted}
  \nucnorm{\PlErrorB} \; \leq \; 3 \nucnorm{\PlErrorA} \, + \, 4
\matsnorm{\Pi_{\SetTwo(\mataspec, \matbspec)}(\ThetaStar)}{1}
\end{equation}
\end{enumerate}
\elems

See Appendix~\ref{AppLemRestricted} for the proof of this claim.
Using Lemma~\ref{LemRestricted}, we can complete the proof of the
theorem.  In particular, from the bound~\eqref{EqnRestricted} and the
RSC assumption, we find that
\begin{align*}
\frac{1}{2 \bigN} \|\Xop(\Delta)\|_2^2 & \geq \kappa(\Xop) \,
\matsnorm{\Delta}{F}^2.
\end{align*}
Using the triangle inequality together with
inequality~\eqref{EqnRestricted}, we obtain
\begin{align*}
\matsnorm{\PlError}{1} & \leq \nucnorm{\PlErrorA} +
\nucnorm{\PlErrorB} \; \leq \; 4 \nucnorm{\PlErrorA} \, + \, 4
\matsnorm{\Pi_{\SetTwo}(\ThetaStar)}{1}.
\end{align*}
From the rank constraint in Lemma~\ref{LemRestricted}(a), we have
$\nucnorm{\PlErrorA} \leq \sqrt{2 \rdim} \matsnorm{\PlErrorA}{F}$.
Putting together the pieces, we find
\begin{align*}
\kappa(\Xop) \, \matsnorm{\PlError}{F}^2 & \leq \max \big \{ 32
\bigreg \: \sqrt{\rdim} \, \matsnorm{\PlError}{F}, \; 16 \: \bigreg \:
\matsnorm{\Pi_{\SetTwo}(\ThetaStar)}{1} \big \},
\end{align*}
which implies that
\begin{align*}
\matsnorm{\PlError}{F} & \leq \max \biggr \{ \frac{32 \bigreg \:
\sqrt{\rdim}}{\kappa(\Xop)}, \; \Big(\frac{16 \: \bigreg \:
\matsnorm{\Pi_{\SetTwo}(\ThetaStar)}{1}}{\kappa(\Xop)} \Big)^{1/2}
\biggr \},
\end{align*}
as claimed.

\subsection{Proof of Corollary~\ref{CorNearLow}}
\label{SecProofCorNearLow}
Let $\tmp = \min \{\lra, \lrb\}$.  In this case, we consider the
singular value decomposition \mbox{$\ThetaStar = U D^T V$,} where $U
\in \real^{\lra \times \tmp}$ and $V \in \real^{\lrb \times \tmp}$ are
orthogonal, and we assume that $D$ is diagonal with the singular
values in non-increasing order $\sigma_1(\ThetaStar) \geq
\sigma_2(\ThetaStar) \geq \ldots \sigma_\tmp(\ThetaStar) \geq 0$.  For
a parameter $\tau > 0$ to be chosen, we let $K = \{i \, \mid \,
\sigma_i(\ThetaStar) > \tau \}$, and let $U^K$ (respectively $V^K$)
denote the $\lra \times |K|$ (respectively the $\lrb \times |K|$)
orthogonal matrix consisting of the first $|K|$ columns of $U$
(respectively $V$).  With this choice, the matrix $\ThetaStar_{\Kcomp}
\mydefn \Pi_{\SetTwo(U^K, V^K)}(\ThetaStar)$ has rank at most $\tmp -
|K|$, with singular values $\{\sigma_i(\ThetaStar), i \in \Kcomp \}$.
Moreover, since $\sigma_i(\ThetaStar) \leq \tau$ for all $i \in
\Kcomp$, we have
\begin{align*}
\nucnorm{\ThetaStar_\Kcomp} & = \tau \sum_{i=|K|+1}^\tmp
\frac{\sigma_i(\ThetaStar)}{\tau} \; \leq \tau \; \sum_{i=|K|+1}^\tmp
\Big(\frac{\sigma_i(\ThetaStar)}{\tau}\Big)^\qpar \; \leq
\tau^{1-\qpar} \radq.
\end{align*}
On the other hand, we also have $\radq \, \geq \; \sum_{i=1}^\tmp
|\sigma_i(\ThetaStar)|^\qpar \; \geq \; |K| \, \tau^\qpar$, which
implies that $|K| \leq \tau^{-\qpar} \, \radq$.  From the general
error bound with $\rdim = |K|$, we obtain
\begin{align*}
\matsnorm{\ThetaHat - \ThetaStar}{F} & \leq \max \biggr \{ \frac{32
\bigreg \: \sqrt{\radq} \, \tau^{-\qpar/2} }{\kappa(\Xop)}, \;
\biggr[\frac{16 \: \bigreg \: \tau^{1-\qpar} \radq }{\kappa(\Xop)}
\biggr]^{1/2} \biggr \},
\end{align*}
Setting $\tau = \bigreg/\kappa$ yields that
\begin{align*}
\matsnorm{\ThetaHat - \ThetaStar}{F} & \leq \max \biggr \{ 
\frac{32
\bigreg^{1-\qpar/2} \: \sqrt{\radq} }{\kappa^{1-\qpar/2}}, \;
\biggr[\frac{16 \: \bigreg^{2-\qpar} \radq }{\kappa^{2-\qpar}}
\biggr]^{1/2} \biggr \} \\
& = 32 \; \sqrt{\radq} \biggr (\frac{\bigreg}{\kappa(\Xop)}
\biggr)^{1-\qpar/2},
\end{align*}
as claimed.
%


\subsection{Proof of Corollary~\ref{CorMultivar}}

For the proof of this corollary, we adopt the following notation.  We
first define the three matrices
\begin{equation}
\label{EqnDefnSpec}
\SpecX = \begin{bmatrix} \Xcov{1}^T \\ \Xcov{2}^T \\ \cdots \\
\Xcov{\numobs}^T \end{bmatrix} \in \real^{\numobs \times \lrb}, \quad
\SpecY = \begin{bmatrix} \Zout{1}^T \\ \Zout{2}^T \\ \cdots \\
\Zout{\numobs}^T \end{bmatrix} \in \real^{\numobs \times \lra}, \quad
\mbox{and} 
\quad
\SpecV = \begin{bmatrix} \Wout{1}^T \\ \Wout{2}^T \\ \cdots \\
\Wout{\numobs}^T \end{bmatrix} \in \real^{\numobs \times \lra}.
\end{equation}
With this notation and using the relation $\bigN = \numobs \lra$, the
SDP objective function~\eqref{EqnSDP} can be written as
$\frac{1}{\lra} \big \{ \frac{1}{2 \numobs} \matsnorm{\SpecY - \SpecX
\Theta}{F}^2 + \regpar \nucnorm{\Theta} \big \}$, where we have
defined $\regpar = \bigreg \, \lra$.

In order to establish the RSC property for this model, some algebra
shows that we need to establish a lower bound on the quantity
\begin{align*}
\frac{1}{2 \numobs} \matsnorm{\SpecX \Delta}{F}^2 & = \frac{1}{2
\numobs} \sum_{j=1}^\pdim \|(\SpecX \Delta)_j\|_2^2 \; \geq \;
\frac{\sigmin(\SpecX^T \SpecX)}{2 \numobs} \: \matsnorm{\Delta}{F}^2,
\end{align*}
where $\sigmin$ denotes the minimum eigenvalue.  The following lemma
follows by adapting known concentration results for random matrices
(see the paper~\cite{Wainwright09a} for details):
\blems
\label{LemMultivarRSC}
Let $\SpecX \in \real^{\numobs \times \pdim}$ be a random matrix with
i.i.d. rows sampled from a $\pdim$-variate $N(0, \Sigma)$
distribution.  Then for $\numobs \geq \pdim$, we have
\begin{equation*}
\mprob \biggr[ \sigmin \Big ( \sampcov \Big ) \geq \frac{\sigmin(\Sigma)}{9}, \;
\sigmax \Big ( \sampcov \Big ) \leq 9 \sigmax(\Sigma) \biggr] \; \geq \; 1 - 4
\exp(-\numobs/2).
\end{equation*}
\elems

\noindent As a consequence, we have $\frac{\sigmin(\SpecX^T \SpecX)}{2
\numobs} \geq \frac{\sigmin(\Sigma)}{18}$ with probability at least
$1- 4 \exp(-\numobs)$ for all $\numobs \geq \pdim$, which establishes
that the RSC property holds with $\kappa(\Xop) = \frac{1}{20}
\sigmin(\Sigma)$.

Next we need upper bound the quantity $\matsnorm{\Xopadj(\mybigW)}{2}$
for this model, so as to verify that the stated choice for $\bigreg$
is valid.  Following some algebra, we find that
\begin{align*}
\frac{1}{\numobs} \, \matsnorm{\Xopadj(\mybigW)}{\myop} & =
\frac{1}{\numobs} \matsnorm{\SpecX^T \SpecV}{\myop}.
\end{align*}
The following lemma is proved in Appendix~\ref{AppLemMultivarNoise}:
\blems
\label{LemMultivarNoise}
There are constants $c_i > 0$ such that
\begin{align}
\mprob \biggr[ \big|\frac{1}{\numobs} \matsnorm{\SpecX^T
\SpecV}{\myop} \big| \geq 5 \nvar \sqrt{\opnorm{\Sigma}}
\sqrt{\frac{\lra + \lrb}{\numobs}} \biggr] & \leq c_1 \exp(-c_2 (\lra
+ \lrb)).
\end{align}
\elems

Using these two lemmas, we can complete the proof of
Corollary~\ref{CorMultivar}.  First, recalling the scaling $\bigN =
\lra \numobs$, we see that Lemma~\ref{LemMultivarNoise} implies that
the choice $\bigreg = 
10 \nvar \sqrt{\opnorm{\Sigma}}
\sqrt{\frac{\lra + \lrb}{\numobs}}$ satisfies the conditions of
Corollary~\ref{CorNearLow} with high probability.
Lemma~\ref{LemMultivarRSC} shows that the RSC property holds with
$\kappa(\Xop) = \sigmin(\Sigma)/20$, again with high probability.
Consequently, Corollary~\ref{CorNearLow} implies that
\begin{align*}
\matsnorm{\ThetaHat - \ThetaStar}{F}^2 & \leq 32^2 \radq \, \biggr (
10 \nvar \sqrt{\opnorm{\Sigma}} \sqrt{\frac{\lra + \lrb}{\numobs}} \;
\; \frac{20}{\sigmin(\Sigma)} \biggr)^{2-\qpar} \\
& = c_1 \; \biggr(\frac{\nvar^2 \opnorm{\Sigma}}{\sigmin^2(\Sigma)}
\biggr)^{1-\qpar/2} \; \radq \; \biggr(\frac{\lra + \lrb}{\numobs}
\biggr)^{1-\qpar/2}
\end{align*}
with probability greater than $1-c_2 \exp(-c_3 (\lra + \lrb))$, as
claimed.

\subsection{Proof of Corollary~\ref{CorAutoregressive}}

For the proof of this corollary, we adopt the notation
\begin{equation*}
\SpecX = \begin{bmatrix} \Zsam{1}^T \\ \Zsam{2}^T \\ \cdots
  \\ \Zsam{\numobs}^T \end{bmatrix} \in \real^{\numobs \times \pdim},
\quad \mbox{and} \quad \SpecY = \begin{bmatrix} \Zsam{2}^T
  \\ \Zsam{2}^T \\ \cdots \\ \Zsam{\numobs+1}^T \end{bmatrix} \in
\real^{\numobs \times \pdim}.
\end{equation*}
Finally, we let $\SpecV \in \real^{\numobs \times \pdim}$ be a matrix
with i.i.d. $N(0, \nvar^2)$ elements, corresponding to the
innovations noise driving the AR process.  With this notation and
using the relation $\bigN = \numobs \pdim$, the SDP objective
function~\eqref{EqnSDP} can be written as $\frac{1}{\pdim} \big \{
\frac{1}{2 \numobs} \matsnorm{\SpecY - \SpecX \Theta}{F}^2 + \regpar
\nucnorm{\Theta} \big \}$, where we have defined $\regpar = \bigreg \,
\pdim$.  At a high level, the proof of this corollary is similar to
that of Corollary~\ref{CorMultivar}, in that we use random matrix
theory to establish the required RSC property, and to justify the
choice of $\regpar$, or equivalently $\bigreg$.  However, it is
considerably more challenging, due to the dependence in the rows of
the random matrices, and the cross-dependence between the two matrices
$\SpecX$ and $\SpecV$ (which were independent in setting of
multivariate regression). \\

The following lemma provides the lower bound needed to establish
RSC for the autoregressive model:
\blems
\label{LemAutoregressiveRSC}
The eigenspectrum of the matrix $\SpecX^T \SpecX/\numobs$ is
well-controlled in terms of the stationary covariance matrix: in
particular, as long as $\numobs > c_3 \pdim$, we have
\begin{equation}
\label{EqnARSig}
  \sigmax \Big (\Sampcov \Big ) \; \stackrel{(a)}{\leq} \; \frac{24
  \,\sigmax(\Sigma)}{1-\maxop}, \quad \mbox{ and } \quad
  \sigmin \Big ( \Sampcov \Big ) \; \stackrel{(b)}{\geq} \;
  \frac{\sigmin(\Sigma)}{4},
\end{equation}
both with probability greater than $1 - 2 c_1 \exp(-c_2 \, \pdim)$.
\elems
\noindent Thus, from the bound~\eqref{EqnARSig}(b), we see with the
high probability, the RSC property holds with $\kappa(\Xop) =
\sigmin(\Sigma)/4$ as long as $\numobs > c_3 \pdim$.

\vspace*{.1in}

As before, in order to verify the choice of $\bigreg$, we need to
control the quantity \mbox{$\frac{1}{\numobs} \matsnorm{\SpecX^T
    \SpecV}{\myop}$.}  The following inequality, proved in
Appendix~\ref{AppLemAutoregressiveNoise}, yields a suitable upper
bound:
\blems
\label{LemAutoregressiveNoise}
There exist constants $c_i > 0$, independent of $\numobs, \pdim,
\Sigma$ etc.  such that
\begin{align}
\label{EqnAutoregressiveNoise}
\mprob \big[\frac{1}{\numobs} \opnorm {\SpecX^T \SpecV} \geq \frac{40
\, \opnorm{\Sigma}}{1-\maxop} \sqrt{\frac{\pdim}{\numobs}} \big] &
\leq c_2 \exp(-c_3 \pdim).
\end{align}
\elems

\noindent From Lemma~\ref{LemAutoregressiveNoise}, we see that it
suffices to choose $\bigreg = \frac{80 \, \opnorm{\Sigma}}{1-\maxop}
\sqrt{\frac{\pdim}{\numobs}}$.  With this choice,
Corollary~\ref{CorNearLow} of Theorem~\ref{ThmMain} yields that
\begin{align*}
\matsnorm{\Theta - \ThetaStar}{F}^2 & \leq c_1 \: \radq \; \biggr[
\frac{\sigmax(\Sigma)}{\sigmin(\Sigma) \, (1 - \maxop)}
\biggr]^{2-\qpar} \, \Big(\frac{\pdim}{\numobs} \Big)^{1-\qpar/2}
\end{align*}
with probability greater than $1- c_2 \exp(-c_3 \pdim)$, as claimed.

\subsection{Proof of Corollary~\ref{CorCompressed}}
\label{SecProofCorCompressed}

Recall that for this model, the observations are of the form
$\mysmally{i} = \tracer{\mysmallX{i}}{\ThetaStar} + \mysmallw{i}$,
where $\ThetaStar \in \real^{\lra \times \lrb}$ is the unknown matrix,
and $\{ \mysmallw{i}\}_{i=1}^\bigN$ is an associated noise sequence.

Let us now show how Proposition~\ref{LemCompressedRSC} implies the RSC
property with an appropriate tolerance parameter.  In particular, let
us define $\delta^2 \; \defn \; \radq \, \big [
\sqrt{\frac{\lra}{\bigN}} + \sqrt{\frac{\lrb}{\bigN}} \big ]^{2 -
\qpar}$, so that if we have the inequality $\matsnorm{\Delta}{F} \;
\leq \; \delta$, the result of Corollary~\ref{CorCompressed} follows
immediately.  Therefore, we may take $\matsnorm{\Delta}{F}^2 \, \geq
\, \delta$. Now recall from Lemma~\ref{LemRestricted} that the error
$\Delta$ satisfies the bound~\eqref{EqnRestricted}.  Combining these
facts, we are guaranteed that $\Error \in \KeySet(\rdim; \delta)$,
where the set $\KeySet$ was previously defined~\eqref{eq:restrict},
and it is sufficient to establish the RSC property over this set.

Observe that the bound~\eqref{EqnCompressedRSC} implies that for any
$\Delta \in \KeySet$,
\begin{align}
  \label{eq:rschelper}
  \frac{\| \Xop(\Error)\|_2}{\sqrt{\bigN}} & \geq \conone
  \matsnorm{\Error}{F} - \contwo \biggr ( \sqrt{\frac{\lra}{\bigN}} \,
  + \, \sqrt{\frac{\lrb}{\bigN}} \biggr ) \nucnorm{\Error}.
\end{align}
Following the arguments used in the proofs of Theorem~\ref{ThmMain}
and Corollary~\ref{CorNearLow}, we find that
\begin{align}
\label{EqnTwoHelper}
\matsnorm{\PlError}{1} \; \leq \; 4 \nucnorm{\PlErrorA} \, + \, 4
\matsnorm{\Pi_{\SetTwo}(\ThetaStar)}{1} \; \leq \; 4 \, \sqrt{2 \radq
\tau^{-\qpar}} \, \matsnorm{\PlErrorA}{F} \, + \, 4 \, \radq
\tau^{1-\qpar},
\end{align}
where $\tau > 0$ is a parameter to be chosen.  We now set $\tau \, =
\, \big (\sqrt{\lra} + \sqrt{\lrb} \, \big )/\sqrt{\bigN}$, and
substitute the resulting bound~\eqref{EqnTwoHelper} into
equation~\eqref{eq:rschelper}, thereby obtaining
\begin{align*}
  \frac{\| \Xop(\Error)\|_2}{\sqrt{\bigN}} & \geq \conone
  \matsnorm{\Error}{F} - \contwo \sqrt{32 \, \radq} \,
  \tau^{1-\qpar/2} \, \frobnorm{\PlError} \, - \, 4 \, \radq \,
  \tau^{2-\qpar}\\
  & \geq \conone \matsnorm{\Error}{F} - \contwo
  \sqrt{32} \, \delta \, \frobnorm{\PlError} \, -
  \, 4 \, \delta \, \frobnorm{\PlError}.
\end{align*}
If we choose $\bigN > 200 \, \radq^{(1-\qpar/2)} \lra \, \lrb$, then
we are guaranteed that $\frac{1}{4} - (4 + \sqrt{32}) \delta \geq
\frac{1}{8}$, which shows that the RSC property holds with
$\kappa(\Xop) = 1/8$. \\

The next step is to control the quantity
$\|\Xopadj(\mybigW)\|_2/\bigN$, required for specifying a suitable
choice of $\bigreg$.
\blems
\label{LemCompressedNoise}
If $\|\mybigW\|_2 \leq 2 \nvar \sqrt{\bigN}$, then
\begin{align}
\mprob \Big[ \frac{\|\Xopadj(\mybigW)\|_2}{\bigN} \geq 4 \nvar
\biggr(\sqrt{\frac{\lra}{\bigN}} + \sqrt{\frac{\lrb}{\bigN}}\biggr)
\Big] & \leq c_1 \exp(-c_2 (\lra + \lrb)).
\end{align}
\elems
\begin{proof}
By definition of the adjoint operator, we have $\frac{1}{\bigN}
\Xopadj(\mybigW) = \frac{1}{\bigN} \sum_{i=1}^\bigN \mysmallw{i}
\mysmallX{i}$.  Since the observation matrices
$\{\mysmallX{i}\}_{i=1}^\bigN$ are i.i.d. Gaussian, if the sequence
$\{\mysmallw{i}\}_{i=1}^\bigN$ is viewed as fixed (by conditioning as
needed), then the random matrix $Z \mydefn \frac{1}{\bigN}
\sum_{i=1}^\bigN \mysmallw{i} \mysmallX{i}$ has zero-mean
i.i.d. Gaussian entries with variance $\frac{\|\mybigW\|^2}{\bigN^2}$.
Since $Z \in \real^{\lra \times \lrb}$, known results in random matrix
theory~\cite{DavSza01} imply that
\begin{align*}
\mprob \biggr[\opnorm{Z} \geq 2 \frac{\|\mybigW\|_2}{\sqrt{\bigN}} \,
\Big( \sqrt{\frac{\lra}{\bigN}} + \sqrt{\frac{\lrb}{\bigN}} \Big)
\biggr] & \leq 2 \exp(-c_2 (\lra + \lrb)),
\end{align*}
as claimed.
\end{proof}


\subsection{Proof of Corollary~\ref{CorRecht}}

This corollary follows from a combination of
Proposition~\ref{LemCompressedRSC} and Lemma~\ref{LemRestricted}.  Let
$\ThetaHat$ be an optimal solution to the
SDP~\eqref{EqnBasisPursuitSDP}, and let $\Delta = \ThetaHat -
\ThetaStar$ be the error.  Since $\ThetaHat$ is optimal and
$\ThetaStar$ is feasible for the SDP, we have $\nucnorm{\ThetaHat} =
\nucnorm{\ThetaStar +\Delta} \leq \nucnorm{\ThetaStar}$.  Using the
decomposition $\PlError = \PlErrorA + \PlErrorB$ from
Lemma~\ref{LemRestricted} and applying triangle inequality, we have
$\nucnorm{\ThetaStar + \PlErrorA + \PlErrorB} \geq \nucnorm{\ThetaStar
+ \PlErrorB} - \nucnorm{\PlErrorA}$.  From the properties of the
decomposition in Lemma~\ref{LemRestricted} (see
Appendix~\ref{AppLemRestricted}), we find that
\begin{align*}
\nucnorm{\ThetaHat} \; = \; \nucnorm{\ThetaStar + \PlErrorA +
\PlErrorB} & \geq \nucnorm{\ThetaStar} + \nucnorm{\PlErrorB} -
\nucnorm{\PlErrorA}.
\end{align*}
Combining the pieces yields that $\nucnorm{\PlErrorB} \leq
\nucnorm{\PlErrorA}$, and hence $\nucnorm{\PlError} \leq 2
\nucnorm{\PlErrorA}$.  By Lemma~\ref{LemRestricted}(a), the rank of
$\PlErrorA$ is at most $2 \rdim$, so that we obtain
$\nucnorm{\PlError} \leq 2 \sqrt{2 \rdim} \matsnorm{\PlError}{F} \leq
4 \rdim \matsnorm{\PlError}{F}$.

Note that $\Xop(\Delta) = 0$, since both $\ThetaHat$ and $\ThetaStar$
agree with the observations.  Consequently, from
Proposition~\ref{LemCompressedRSC}, we have that
\begin{align*}
0 \; = \; \frac{\| \Xop(\Delta)\|_2}{\sqrt{\bigN}} & \geq \conone
\matsnorm{\Delta}{F} - \contwo \BigMess \nucnorm{\Delta} \\
& \geq \matsnorm{\Delta}{F} \Biggr( \conone - 4 \sqrt{\frac{\rdim
\lra}{\bigN}} + 4 \sqrt{\frac{\rdim \lrb}{\bigN}} \Biggr) \\
& \geq \matsnorm{\Delta}{F}/20,
\end{align*}
where the final inequality follows from the assumption that $\bigN >
40^2 \rdim (\lra + \lrb)$.  We have thus shown that $\Delta = 0$,
which implies that $\ThetaHat = \ThetaStar$ as claimed.


\section{Experimental results}
\label{SecSimulations}

In this section, we report the results of various simulations that
demonstrate the close agreement between the scaling predicted by our
theory, and the actual behavior of the SDP-based
$M$-estimator~\eqref{EqnSDP} in practice.  In all cases, we solved the
convex program~\eqref{EqnSDP} by using our own implementation in
MATLAB of an accelerated gradient descent method which adapts a
non-smooth convex optimization procedure~\cite{Nesterov07} to the
nuclear-norm~\cite{Yi09}.  We chose the regularization parameter
$\bigreg$ in the manner suggested by our theoretical results; in doing
so, we assumed knowledge of quantities such as the noise variance
$\nvar^2$.  (In practice, one would have to estimate such quantities
from the data using standard methods.)

We report simulation results for three of the running examples
discussed in this paper: low-rank multivariate regression, estimation
in vector autoregressive processes, and matrix recovery from random
projections (compressed sensing).  In each case, we solved instances
of the SDP for a square matrix $\ThetaStar \in \real^{\lrb \times
  \lrb}$, where $\lrb \in \{40, 80, 160 \}$ for the first two
examples, and $\lrb \in \{20, 40, 80\}$ for the compressed sensing
example.  In all cases, we considered the case of exact low rank
constraints, with $\rank(\ThetaStar) = \rdim = 10$, and we generated
$\ThetaStar$ by choosing the subspaces of its left and right singular
vectors uniformly at random from the Grassman manifold.  The
observation or innovations noise had variance $\nvar^2 = 1$ in each
case.  The VAR process was generated by first solving for the
covariance matrix $\Sigma$ using the MATLAB function dylap and then
generating a sample path.  For each setting of $(\rdim, \pdim)$, we
solved the SDP for a range of sample sizes $\bigN$.

\begin{figure}[h]
  \begin{center}
    \begin{tabular}{ccc}
     \widgraph{.45\textwidth}{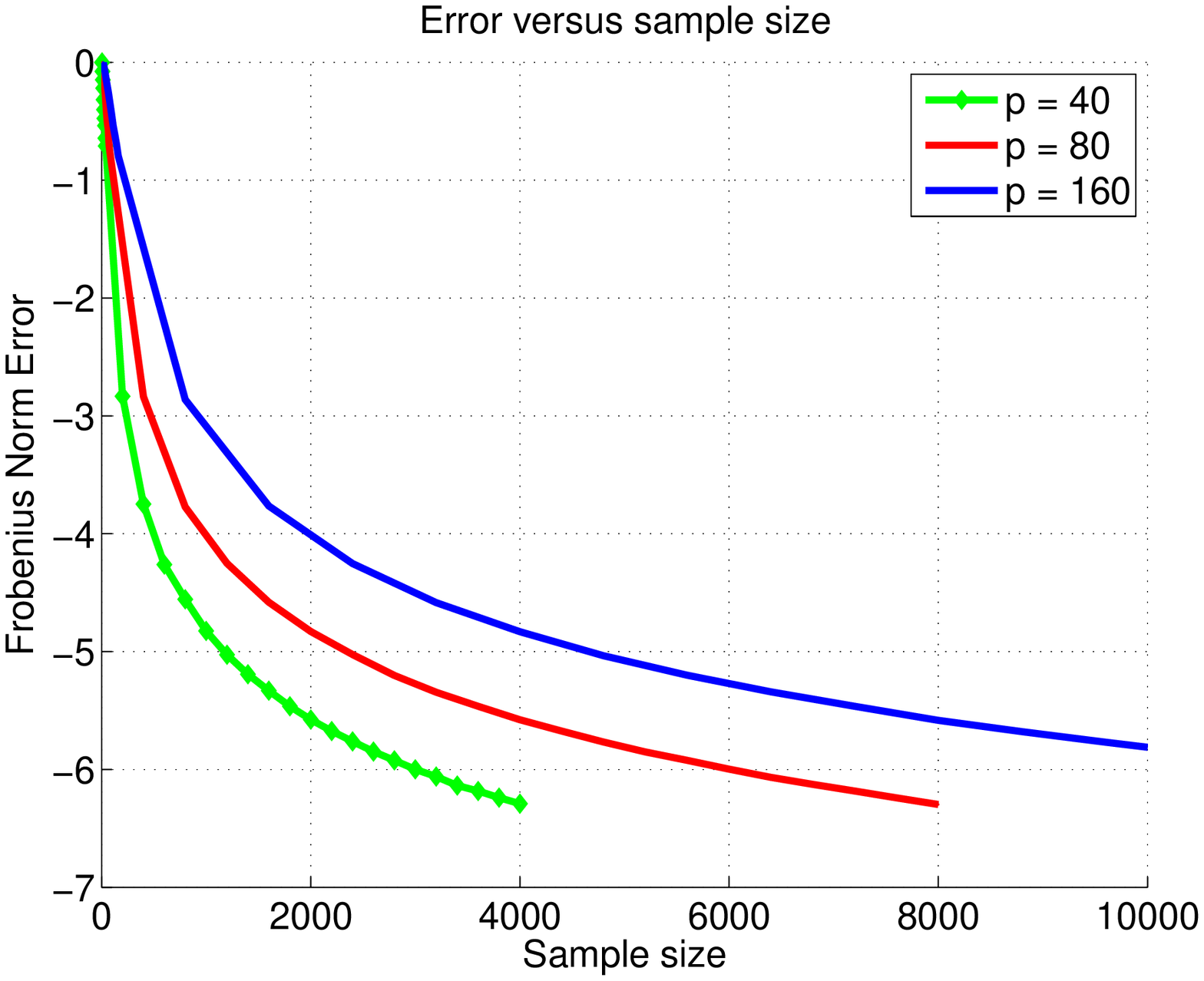} & & 
     \widgraph{.45\textwidth}{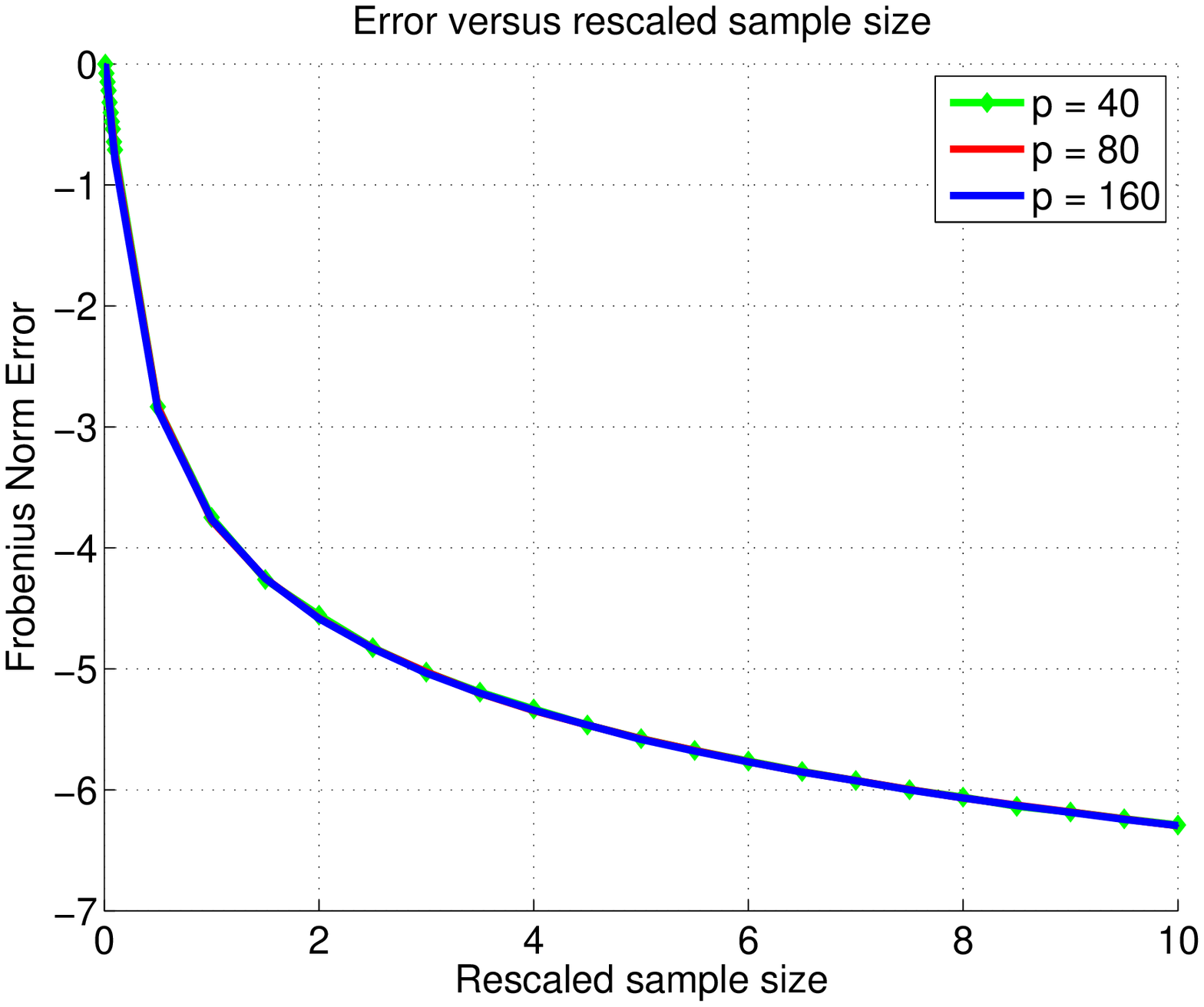} \\
      (a) & & (b)
    \end{tabular}
  \end{center}
  \caption{Results of applying the SDP~\eqref{EqnSDP} with nuclear
norm regularization to the problem of low-rank multivariate
regression.  (a) Plots of the Frobenius error $\matsnorm{\ThetaHat -
\ThetaStar}{F}$ on a logarithmic scale versus the sample size $\bigN$
for three different matrix sizes $\lrb \in \{40, 80, 160\}$, all with
rank $\rdim = 10$.  (b) Plots of the same Frobenius error versus the
rescaled sample size $\bigN/(\rdim \pdim)$.  Consistent with theory,
all three plots are now extremely well-aligned.}
\label{FigMultivarResults}
\end{figure}

Figure~\ref{FigMultivarResults} shows results for a multivariate
regression model with the covariates chosen randomly from a $N(0, I)$
distribution.  Panel (a) plots the Frobenius error
$\matsnorm{\ThetaHat - \ThetaStar}{F}$ on a logarithmic scale versus
the sample size $\bigN$ for three different matrix sizes, $\lrb \in
\{40, 80, 160 \}$.  Naturally, in each case, the error decays to zero
as $\bigN$ increases, but larger matrices require larger sample sizes,
as reflected by the rightward shift of the curves as $\lrb$ is
increased.  Panel (b) of Figure~\ref{FigMultivarResults} shows the
exact same set of simulation results, but now with the Frobenius error
plotted versus the rescaled sample size $\widetilde{\bigN} \mydefn
\bigN/(\rdim \lrb)$.  As predicted by Corollary~\ref{CorMultivar}, the
error plots now are all aligned with one another; the degree of
alignment in this particular case is so close that the three plots are
now indistinguishable.  (The blue curve is the only one visible since
it was plotted last by our routine.)  Consequently,
Figure~\ref{FigMultivarResults} shows that $\bigN/(\rdim \lrb)$ acts
as the effective sample size in this high-dimensional setting.

\begin{figure}[h]
  \begin{center}
\begin{tabular}{ccc}
     \widgraph{.45\textwidth}{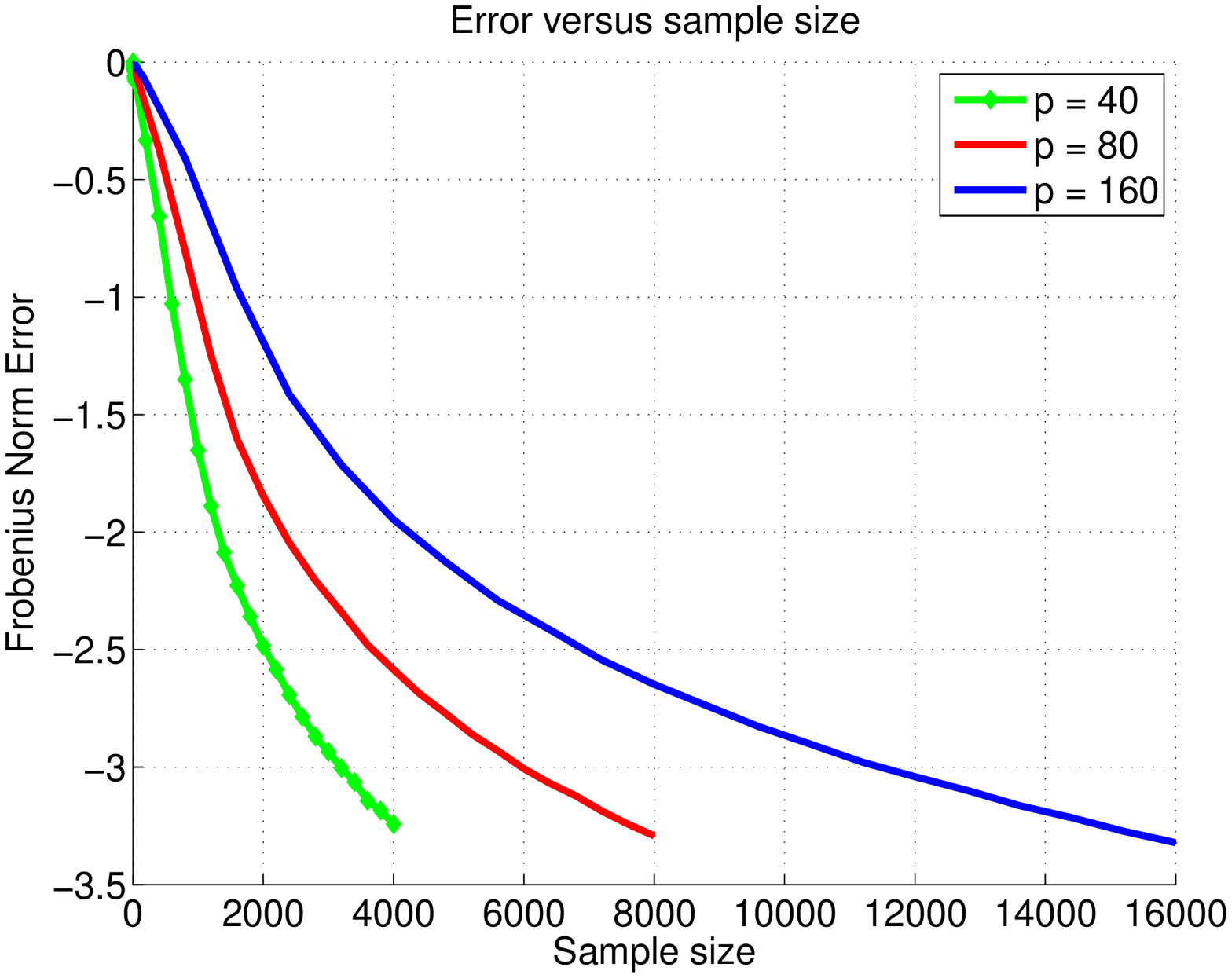} & & 
     \widgraph{.45\textwidth}{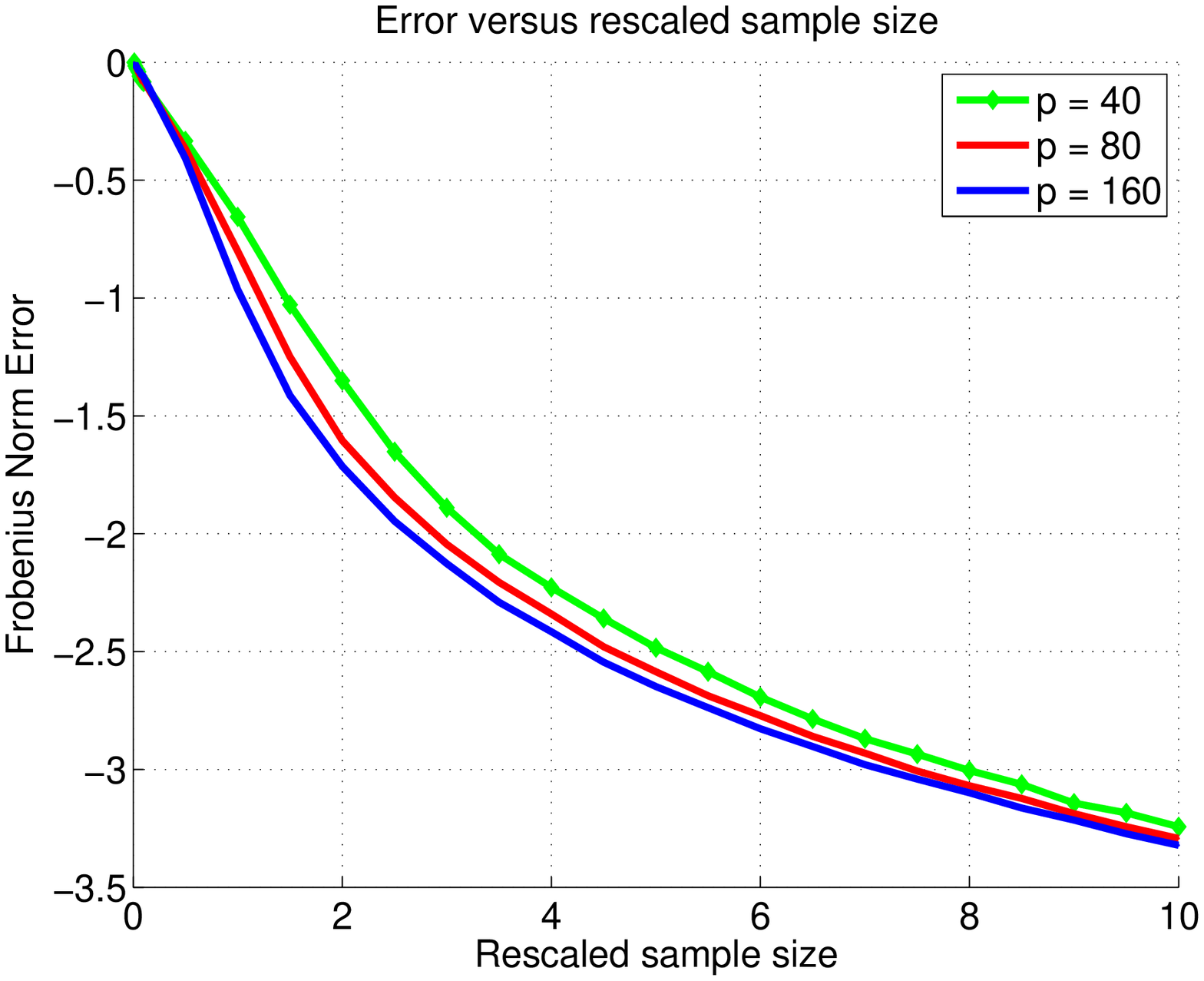} \\
(a) & & (b)
\end{tabular}
\end{center}
\caption{Results of applying the SDP~\eqref{EqnSDP} with nuclear norm
regularization to estimating the system matrix of a vector
autoregressive process.  (a) Plots of the Frobenius error
$\matsnorm{\ThetaHat - \ThetaStar}{F}$ on a logarithmic scale versus
the sample size $\bigN$ for three different matrix sizes $\lrb \in
\{40, 80, 160\}$, all with rank $\rdim = 10$.  (b) Plots of the same
Frobenius error versus the rescaled sample size $\bigN/(\rdim \pdim)$.
Consistent with theory, all three plots are now reasonably
well-aligned.}
\label{FigAutoregressiveResults}
\end{figure}

Figure~\ref{FigAutoregressiveResults} shows similar results for the
autoregressive model discussed in Example~\ref{ExaAutoregressive}.  As
shown in panel (a), the Frobenius error again decays as the sample
size is increased, although problems involving larger matrices are
shifted to the right.  Panel (b) shows the same Frobenius error
plotted versus the rescaled sample size $\bigN/(\rdim \lrb)$; as
predicted by Corollary~\ref{CorAutoregressive}, the errors for
different matrix sizes $\lrb$ are again quite well-aligned.  In this
case, we find (both in our theoretical analysis and experimental
results) that the dependence in the autoregressive process slows down
the rate at which the concentration occurs, so that the results are
not as crisp as the low-rank multivariate setting in
Figure~\ref{FigMultivarResults}.

\begin{figure}[h]
  \begin{center}
    \begin{tabular}{ccc}
     \widgraph{.45\textwidth}{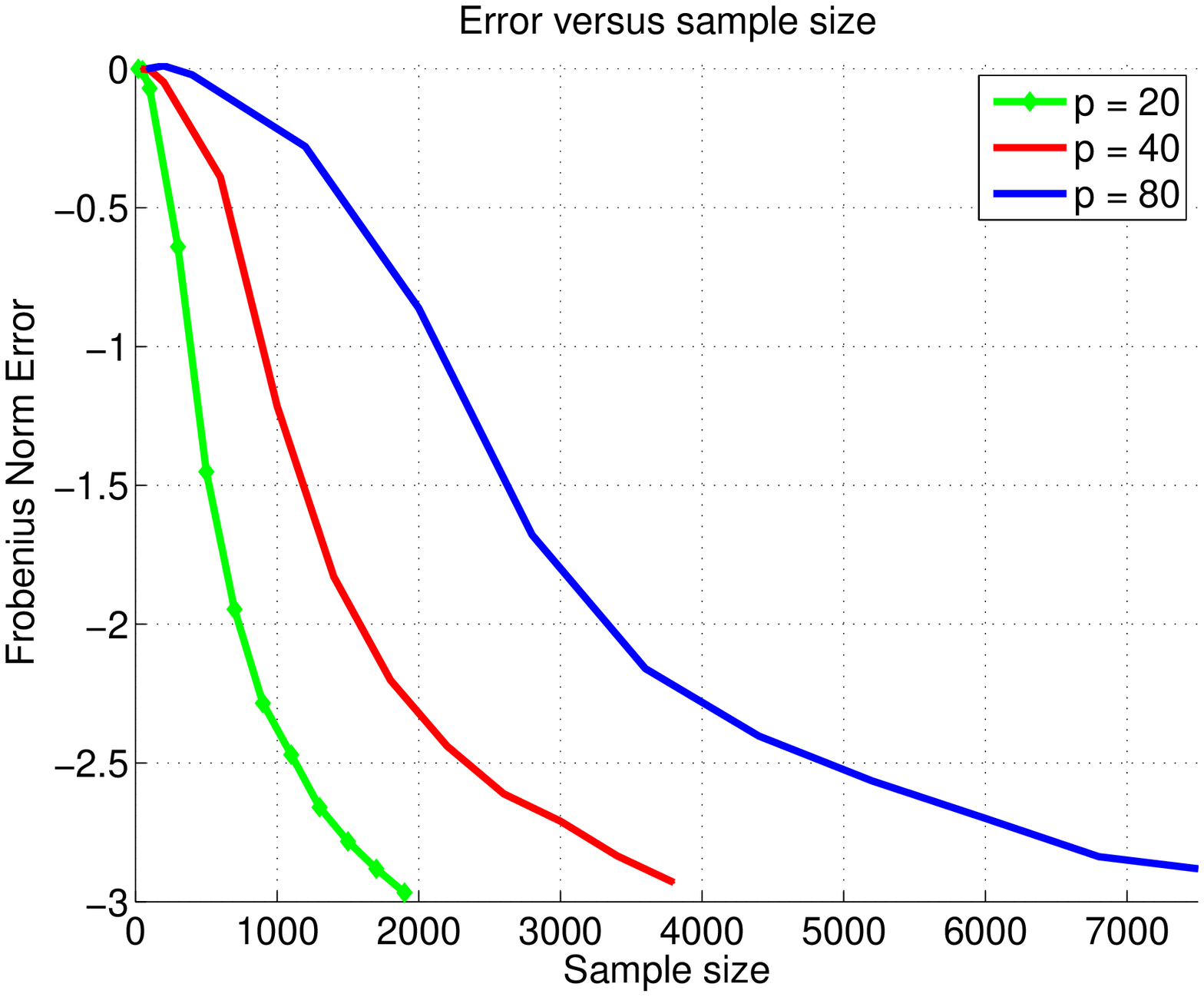} & &
     \widgraph{.45\textwidth}{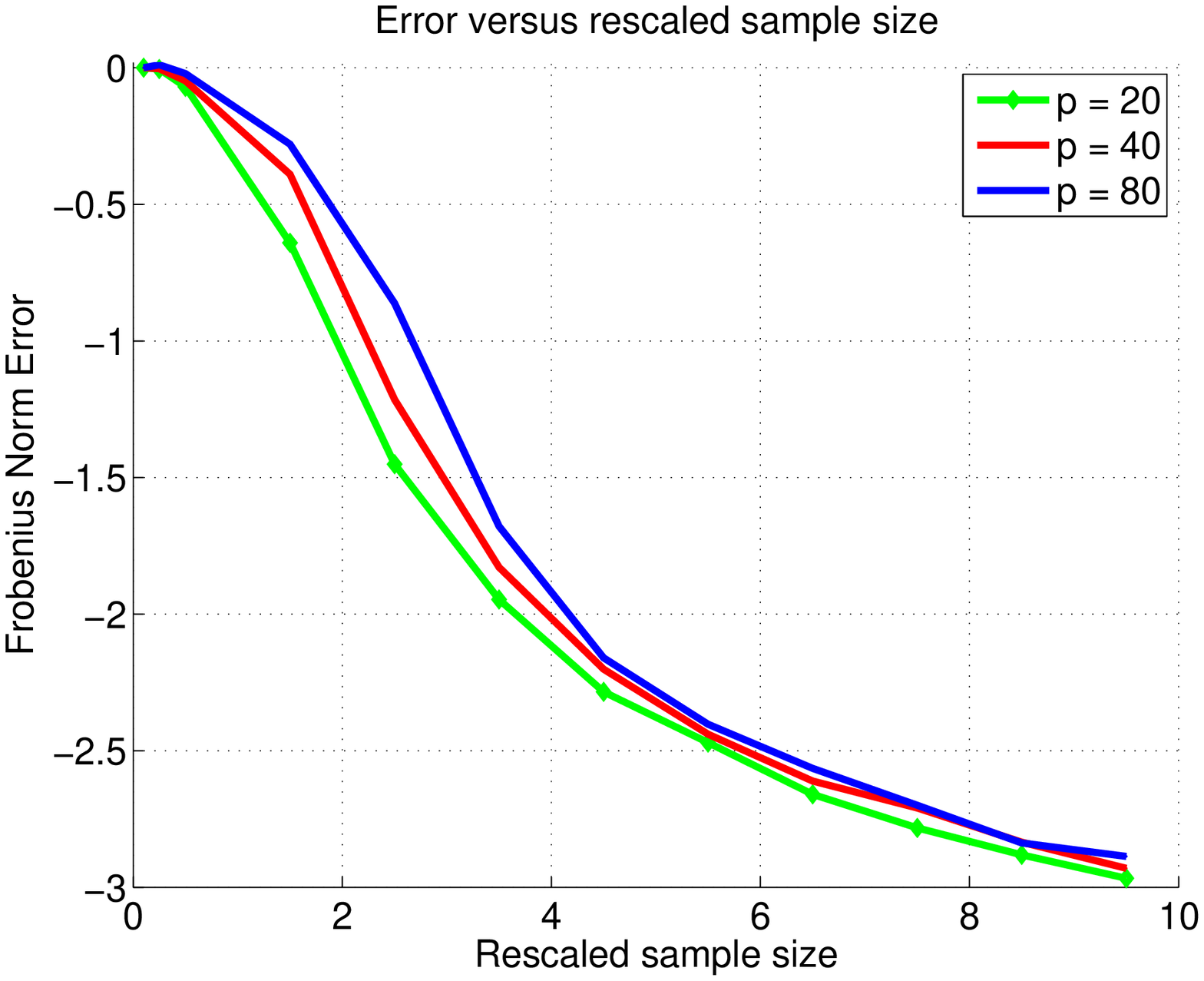} \\
      (a) & & (b)
    \end{tabular}
  \end{center}
\caption{ Results of applying the SDP~\eqref{EqnSDP} with nuclear norm
regularization to recovering a low-rank matrix on the basis of random
projections (compressed sensing model) (a) Plots of the Frobenius
error $\matsnorm{\ThetaHat - \ThetaStar}{F}$ on a logarithmic scale
versus the sample size $\bigN$ for three different matrix sizes $\lrb
\in \{20, 40, 80\}$, all with rank $\rdim = 10$.  (b) Plots of the
same Frobenius error versus the rescaled sample size $\bigN/(\rdim
\pdim)$.  Consistent with theory, all three plots are now reasonably
well-aligned.}
\label{FigCompressedResults}
\end{figure}

Finally, Figure~\ref{FigCompressedResults} presents the same set of
results for the compressed sensing observation model discussed in
Example~\ref{ExaCompressed}.  Even though the observation matrices
$X_i$ here are qualitatively different (in comparison to the
multivariate regression and autoregressive examples), we again see the
``stacking'' phenomenon of the curves when plotted versus the rescaled
sample size $\bigN / \rdim \lrb$, as predicted by
Corollary~\ref{CorCompressed}.

\section{Discussion}

In this paper, we have analyzed the nuclear norm relaxation for a
general class of noisy observation models, and obtained non-asymptotic
error bounds on the Frobenius norm that hold under high-dimensional
scaling.  In contrast to most past work, our results are applicable to
both exactly and approximately low-rank matrices.  We stated a main
theorem that provides high-dimensional rates in a fairly general
setting, and then showed how by specializing this result to some
specific model classes---namely, low-rank multivariate regression,
estimation of autoregressive processes, and matrix recovery from
random projections---it yields concrete and readily interpretable
rates.  Lastly, we provided some simulation results that showed
excellent agreement with the predictions from our theory.

This paper has focused on achievable results for low-rank matrix
estimation using a particular polynomial-time method.  It would be
interesting to establish matching lower bounds, showing that the rates
obtained by this estimator are minimax-optimal.  We suspect that this
should be possible, for instance by using the techniques exploited in
Raskutti et al.~\cite{RasWaiYu09} in analyzing minimax rates for
regression over $\ell_\qpar$-balls.

\label{SecDiscussion}


\section*{Acknowledgements}

This work was partially supported by a Sloan Foundation Fellowship,
AFOSR-09NL184 grant, and an NSF-CCF-0545862 CAREER grant to MJW.

\appendix

\section{Proof of Lemma~\ref{LemRestricted}}
\label{AppLemRestricted}

Part (a) of the claim was proved in Recht et al.~\cite{Recht07}; we
simply provide a proof here for completeness.  We write the SVD as
$\ThetaStar = U D V^T$, where $U \in \real^{\lra \times \lra}$ and $V
\in \real^{\lrb \times \lrb}$ are orthogonal matrices, and $D$ is the
matrix formed by the singular values of $\ThetaStar$.  By re-ordering
as needed, we may assume without loss of generality that the first
$\rdim$ columns of $U$ (respectively $V$) correspond to the matrices
$\mataspec$ (respectively $\matbspec$) from the statement.  We then
define the matrix $\DelRotate = U^T \Delta V \in \real^{\lra \times
\lrb}$, and write it in block form as
\begin{align*}
\DelRotate & = \begin{bmatrix} \DelRotate_{11} & \DelRotate_{12} \\
  \DelRotate_{21} & \DelRotate_{22}
\end{bmatrix}, \qquad \mbox{where 
$\DelRotate_{11} \in \real^{\rdim \times \rdim}$, and $\DelRotate_{22}
\in \real^{(\lra - \rdim) \times (\lrb - \rdim)}$.}
\end{align*}                     
We now define the matrices
\begin{align*}
\PlErrorB & = U \begin{bmatrix} 0 & 0\\ 0 & \DelRotate_{22}
  \end{bmatrix} V^T, \qquad \mbox{ and } \PlErrorA = \PlError - \PlErrorB.
\end{align*}
Note that we have
\begin{align*}
\rank(\PlErrorA) & = \rank \begin{bmatrix} \DelRotate_{11} &
\DelRotate_{12} \\ \DelRotate_{21} & 0 \end{bmatrix} \; \leq \; \rank
\begin{bmatrix} \DelRotate_{11} & \DelRotate_{12} \\ 0 & 0
\end{bmatrix} + \rank \begin{bmatrix} \DelRotate_{11} & 0 \\
\DelRotate_{21} & 0 \end{bmatrix} \; \leq 2 \rdim,
\end{align*}
which establishes Lemma~\ref{LemRestricted}(a).  Moreover, we note for
future reference that by construction of $\PlErrorB$, the nuclear
norm satisfies the decomposition
\begin{align}
\label{EqnDecompose}
\nucnorm{\Pi_{\SetOne(\mataspec, \matbspec)}(\ThetaStar) + \PlErrorB}
& = \nucnorm{\Pi_{\SetOne(\mataspec, \matbspec)}(\ThetaStar)} +
\nucnorm{\PlErrorB}.
\end{align}

We now turn to the proof of Lemma~\ref{LemRestricted}(b).  
Recall that the error $\PlError = \ThetaHat - \ThetaStar$
associated with any optimal solution must satisfy the
inequality~\eqref{EqnBasic}, which implies that
\begin{align}
\label{EqnImply}
0 & \leq \frac{1}{\bigN} \inprod{\mybigW}{\Xop(\PlError)} + \bigreg
\big \{ \matsnorm{\ThetaStar}{1} - \matsnorm{\ThetaHat}{1} \big \} \;
\leq \; \opnorm{\frac{1}{\bigN} \Xopadj(\mybigW)} \,
\nucnorm{\PlError} + \bigreg \big \{ \matsnorm{\ThetaStar}{1} -
\matsnorm{\ThetaHat}{1} \big \},
\end{align}
where we have used the bound~\eqref{EqnCafeOne}.

Using the triangle inequality and the relation~\eqref{EqnDecompose}, we have
\begin{align*}
\nucnorm{\ThetaHat} & = \nucnorm{(\Pi_{\SetOne}(\ThetaStar) +
\PlErrorB) + (\Pi_{\SetTwo}(\ThetaStar) + \PlErrorA)} \\
& \geq \nucnorm{(\Pi_{\SetOne}(\ThetaStar) + \PlErrorB)} -
\nucnorm{(\Pi_{\SetTwo}(\ThetaStar) + \PlErrorA)} \\
& \geq \nucnorm{\Pi_{\SetOne}(\ThetaStar)} + \nucnorm{\PlErrorB}  -
\big \{ \nucnorm{(\Pi_{\SetTwo}(\ThetaStar)} + \nucnorm{\PlErrorA} \big \}.
\end{align*}
Consequently, we have
\begin{align*}
\matsnorm{\ThetaStar}{1} - \matsnorm{\ThetaHat}{1} & \leq
\nucnorm{\ThetaStar} - \big \{ \nucnorm{\Pi_{\SetOne}(\ThetaStar)} +
\nucnorm{\PlErrorB} \big \} + \big \{
\nucnorm{(\Pi_{\SetTwo}(\ThetaStar)} + \nucnorm{\PlErrorA} \big\} \\
& = 2 \nucnorm{\Pi_{\SetTwo}(\ThetaStar)} + \nucnorm{\PlErrorA} -
\nucnorm{\PlErrorB}.
\end{align*}
Substituting this inequality into the bound~\eqref{EqnImply}, we obtain
\begin{align*}
0 & \leq \opnorm{\frac{1}{\bigN} \Xopadj(\mybigW)} \,
\nucnorm{\PlError} + \bigreg \big \{ 2
\nucnorm{\Pi_{\SetTwo}(\ThetaStar)} + \nucnorm{\PlErrorA} -
\nucnorm{\PlErrorB} \big \}.
\end{align*}
Finally, since $\opnorm{\frac{1}{\bigN} \Xopadj(\mybigW)} \leq \bigreg/2$
by assumption, we conclude that
\begin{align*}
0 & \leq \bigreg \big \{ 2 \nucnorm{\Pi_{\SetTwo}(\ThetaStar)} +
\frac{3}{2} \nucnorm{\PlErrorA} - \frac{1}{2} \nucnorm{\PlErrorB} \big
\},
\end{align*}
from which the bound~\eqref{EqnRestricted} follows.


\section{Proof of Lemma~\ref{LemMultivarNoise}}
\label{AppLemMultivarNoise}

Let $\Sphere{\pdim} = \{ u \in \real^\pdim \, \mid \, \|u\|_2 = 1 \}$
denote the Euclidean sphere in $\pdim$-dimensions.  The operator norm
of interest has the variational representation
\begin{align*}
\frac{1}{\numobs} \opnorm{\SpecX^T \SpecV} & = \frac{1}{\numobs}
\sup_{u \in \Sphere{\lra}} \; \sup_{v \in \Sphere{\lrb}} v^T \SpecX^T
\SpecV u
\end{align*}
For positive scalars $a$ and $b$, define the (random) quantity
\begin{align*}
\Psi(a,b) & \mydefn  \sup_{u \in a \, \Sphere{\lra}} \sup_{v
\in b \, \Sphere{\lrb}} \inprod{\SpecX v}{\SpecV u}.
\end{align*}
and note that our goal is to upper bound $\Psi(1,1)$.  Note moreover
that $\Psi(a,b) = a \, b \, \Psi(1,1)$, a relation which will be
useful in the analysis.

Let $\Covera = \{u^1, \ldots, u^\covnuma \}$ and $\Coverb = \{ v^1,
\ldots, v^\covnumb \}$ denote $1/4$ coverings of $\Sphere{\lra}$ and
$\Sphere{\lrb}$, respectively.  We now claim that we have the upper
bound
\begin{equation}
\label{EqnInterClaim}
  \Psi(1,1) \; \leq \; 4 \max_{u^a \in \Covera, v^b \in \Coverb}
\inprod{\SpecX v^b}{\SpecV u^a}
\end{equation}
To establish this claim, we note that since the sets $\Covera$ and
$\Coverb$ are $1/4$-covers, for any pair $(u,v) \in \Sphere{\pdim}
\times \Sphere{\pdim}$, there exists a pair $(u^a, v^b) \in \Covera
\times \Coverb$ such that $u=u^a + \Delta u$ and $v=v^b + \Delta v$,
with $\max \{ \|\Delta u\|_2, \|\Delta v\|_2 \} \leq 1/4$.
Consequently, we can write
\begin{align}
\label{EqnFourWay}
\inprod{\SpecX v}{\SpecV u} & = \inprod{\SpecX v^b}{\SpecV u^a} +
\inprod{\SpecX v^b}{\SpecV \Delta u} + \inprod{\SpecX \Delta v}{\SpecV
u^a} + \inprod{\SpecX \Delta v}{\SpecV \Delta u}.
\end{align}
By construction, we have the bound $|\inprod{\SpecX v^b}{\SpecV \Delta
u}| \leq \Psi(1, 1/4) = \frac{1}{4} \Psi(1,1)$, and similarly
$|\inprod{\SpecX \Delta v}{\SpecV u^a}| \leq \frac{1}{4} \Psi(1,1)$ as
well as $|\inprod{\SpecX \Delta v}{\SpecV \Delta u}| \leq \frac{1}{16}
\Psi(1,1)$.  Substituting these bounds into the
decomposition~\eqref{EqnFourWay} and taking suprema over the left and
right-hand sides, we conclude that
\begin{align*}
\Psi(1,1) & \leq \max_{u^a \in \Covera, v^b \in \Coverb}
\inprod{\SpecX v^b}{\SpecV u^a}  + \frac{9}{16} \Psi(1,1),
\end{align*}
from which the bound~\eqref{EqnInterClaim} follows.

We now apply the union bound to control the discrete maximum.  It is
known (e.g.,~\cite{LedTal91,Matousek}) that there exists a $1/4$
covering of $\Sphere{\lra}$ and $\Sphere{\lrb}$ with at most $\covnuma
\leq 8^\lra$ and $\covnumb \leq 8^\lrb$ elements respectively.
Consequently, we have
\begin{align}
\label{EqnUnion}
\mprob \big[|\Psi(1,1)| \geq 4 \delta \, \numobs \big] & \leq 8^{\lra
+ \lrb} \max_{u^a, v^b} \; \; \mprob \left [ \frac{ |\inprod{\SpecX
v^b}{\SpecV u^a}|}{\numobs} \geq \delta \right ].
\end{align}
It remains to obtain a good bound on the quantity $\frac{1}{\numobs}
\inprod{\SpecX v}{\SpecV u} = \frac{1}{\numobs} \sum_{i=1}^\numobs
\inprod{v}{\SpecX_i} \inprod{u}{\SpecV_i}$, where $(u,v) \in
\Sphere{\lra} \times \Sphere{\lrb}$ are arbitrary but fixed.  Since
$\SpecV_i \in \real^\lra$ has i.i.d. $N(0, \nvar^2)$ elements and $u$
is fixed, we have $Z_i \mydefn \inprod{u}{\SpecV_i} \sim N(0,
\nvar^2)$ for each $i = 1, \ldots, \numobs$.  These variables are
independent of one another, and of the random matrix $\SpecX$.
Therefore, conditioned on $\SpecX$, the sum $Z \mydefn
\frac{1}{\numobs} \sum_{i=1}^\numobs \inprod{v}{\SpecX_i}
\inprod{u}{\SpecV_i}$ is zero-mean Gaussian with variance
\begin{align*}
\alpha^2 & \mydefn \frac{\nvar^2}{\numobs} \left (\frac{1}{\numobs}
\|\SpecX v\|_2^2 \right ) \; \leq \; 
\frac{\nvar^2}{\numobs}
\opnorm{\SpecX^T \SpecX/\numobs}.
\end{align*}
Define the event $\Tail = \{\alpha^2 \leq \frac{9 \nvar^2
\opnorm{\Sigma}}{\numobs} \}$.  Using Lemma~\ref{LemMultivarRSC}, we
have $\opnorm{\SpecX^T \SpecX/\numobs} \leq 9 \sigmax(\Sigma)$ with
probability at least $1 - 2 \exp(-\numobs/2)$, which implies that
$\mprob[\Tail^c] \leq 2 \exp(-\numobs/2)$.  Therefore, conditioning on
the event $\Tail$ and its complement $\Tail^c$, we obtain
\begin{align*}
\mprob[|Z| \geq t] & \leq \mprob \big[ |Z| \geq t \mid \Tail \big] +
\mprob[\Tail^c] \\
& \leq \exp \left (-\numobs \frac{t^2}{2 \nvar^2 \,(4 +
\opnorm{\Sigma})} \right ) + 2 \exp(-\numobs/2).
\end{align*}
Combining this tail bound with the upper bound~\eqref{EqnUnion}, we
have
\begin{align*}
\mprob \big[|\psi(1,1)| \geq 4 \delta \, \numobs \big] & \leq 8^{\lra
+ \lrb} \left \{ \exp \left (-\numobs \frac{t^2}{18 \nvar^2
\opnorm{\Sigma}} \right) + 2 \exp(-\numobs/2) \right \}.
\end{align*}
Setting $t^2 = 20 \nvar^2 \opnorm{\Sigma} \frac{\lra +
\lrb}{\numobs}$, this probability vanishes as long as $\numobs > 16
(\lra + \lrb)$.


\section{Technical details for Corollary~\ref{CorAutoregressive}}

In this appendix, we collect the proofs of
Lemmas~\ref{LemAutoregressiveRSC} and~\ref{LemAutoregressiveNoise}.

\subsection{Proof of Lemma~\ref{LemAutoregressiveRSC}}
\label{AppLemAutoregressiveRSC}
Recalling that $\Sphere{\pdim}$ denotes the unit-norm Euclidean sphere
in $\pdim$-dimensions, we first observe that $\opnorm{X} = \sup_{u \in
\Sphere{\pdim}} \|X u\|_2$.  Our next step is to reduce the supremum
to a maximization over a finite set, using a standard covering
argument.  Let \mbox{$\Covera = \{u^1, \ldots, u^\covnuma \}$} denote
a $1/2$-cover of it.  By definition, for any $u \in \Sphere{\pdim}$,
there is some $u^a \in \Covera$ such that $u = u^a + \Delta u$, where
$\|\Delta u\|_2 \leq 1/2$.  Consequently, for any $u \in
\Sphere{\pdim}$, the triangle inequality implies that
\begin{align*}
\|X u\|_2 & \leq  \|X u^a\|_2 + \|X \Delta u\|_2,
\end{align*}
and hence that
$\opnorm{X} \leq \max_{u^a \in \Covera} \|X u^a\|_2 + \frac{1}{2} \opnorm{X}$.
 Re-arranging yields the useful inequality
\begin{align}
\label{EqnDiscreteOne}
\opnorm{X} \leq 2 \max_{u^a \in \Covera} \|X u^a\|_2.
\end{align}

Using inequality~\eqref{EqnDiscreteOne}, we have
\begin{align}
  \mprob \biggr[\frac{1}{\numobs} \opnorm{\SpecX^T \SpecX} > t \biggr]
  & \leq \mprob \biggr [\max_{u^a \in \Covera} \frac{1}{\numobs}
  \sum_{i=1}^\numobs (\inprod{u^a}{\SpecX_i})^2 > \frac{t}{2} \biggr ]
  \nonumber \\
& 
\label{EqnUnionRoma}
\leq 4^\pdim \, \max_{u^a \in \Covera} \, 
\mprob \biggr [
    \frac{1}{\numobs} \sum_{i=1}^\numobs (\inprod{u^a}{\SpecX_i})^2 >
    \frac{t}{2} \biggr ].
\end{align}
where the last inequality follows from the union bound, and the
fact~\cite{LedTal91,Matousek} that there exists a $1/2$-covering of
$\Sphere{\pdim}$ with at most $4^\pdim$ elements.

  In order to complete the proof, we need to obtain a sharp upper
bound on the quantity \mbox{$\mprob \big[ \frac{1}{\numobs}
\sum_{i=1}^\numobs (\inprod{u}{\SpecX_i})^2 > \frac{t}{2} \big]$,}
valid for any fixed $u \in \Sphere{\pdim}$.  Define the random vector
$Y \in \real^\numobs$ with elements $Y_i = \binprod{u}{\SpecX_i}$.
Note that $Y$ is zero mean, and its covariance matrix $R$ has elements
\mbox{$R_{ij} = \Exs[Y_i Y_j] \, = \, u^T \Sigma (\og)^{|j-i|} \, u$.}
In order to bound the spectral norm of $R$, we note that since it is
symmetric, we have \mbox{$\opnorm{R} \leq \max \limits_{i=1, \ldots,
\pdim} \sum_{j=1}^\pdim |R_{ij}|$,} and moreover
\begin{align*}
|R_{ij}| & = |u^T \Sigma (\og)^{|j-i|} \, u| \; \leq \;
(\opnorm{\ThetaStar})^{|j-i|} \op{\Sigma} \; \leq \; \maxop^{|j-i|} \,
\opnorm{\Sigma}.
\end{align*}
Combining the pieces, we obtain
\begin{align}
\label{EqnRBound}
\opnorm{R} & \leq \max_i \sum_{j=1}^\pdim |\maxop|^{|i-j|}
\opnorm{\Sigma} \; \leq \; 2 \opnorm{\Sigma} \sum_{j=0}^\infty
|\maxop|^{j} \; \leq \; \frac{2 \opnorm{\Sigma}}{1-\maxop}.
\end{align}
Moreover, we have $\trace(R)/\numobs = u^T \Sigma u \leq
\opnorm{\Sigma}$.
Applying Lemma~\ref{LemMeta} with $t = 5 \sqrt{\frac{\pdim}{\numobs}}$,
we conclude that
\begin{align*}
  \mprob \biggr[ \, \frac{1}{\numobs} \, \| Y\|_2^2 \; > \;
    \opnorm{\Sigma} + \, 5 \sqrt{\frac{\pdim}{\numobs}} \, \opnorm{R}
    \biggr] \; \leq \; 2 \exp \big ( - 5 \pdim \big ) + 2 \exp -\numobs/2 )..
\end{align*}

Combined with the bound~\eqref{EqnUnionRoma}, we obtain
\begin{equation}
  \opnorm{\sampcov} \, \leq \, \opnorm{\Sigma} \biggr \{ 2 +
  \frac{20}{ (1-\maxop)} \, \sqrt{\frac{\pdim}{\numobs}} \biggr \} \; \leq \;
\frac{24 \opnorm{\Sigma}}{(1-\maxop)},
\end{equation}
with probability at least $1 - c_1 \, \exp(-c_2 \, \pdim)$, which
establishes the upper bound~\eqref{EqnARSig}(a).

\vspace*{.2in}

Turning to the lower bound~\eqref{EqnARSig}(b), we let $\Coverb =
\{v^1, \ldots, v^\covnumb \}$ be an $\epsilon$-cover of
$\Sphere{\pdim}$ for some $\epsilon \in (0,1)$ to be chosen.  Thus,
for any $v \in \real^\pdim$, there exists some $v^b$ such that $v =
v^b + \Delta v$, and $\| \Delta v \|_2 \leq \epsilon$.  Define the
function $\Psi: \real^\pdim \times \real^\pdim \rightarrow \real$ via
$\Psi(u,v) = u^T \Sampcov v$, and note that $\Psi(u,v) = \Psi(v,u)$.
With this notation, we have
\begin{align*}
v^T \Sampcov v = \Psi(v, v) & = \Psi(v^k,v^k) + 2 \Psi(\Delta v,v) +
\Psi(\Delta v,\Delta v) \\
& \geq \Psi(v^k,v^k) + 2 \Psi(\Delta v,v),
\end{align*}
since $\Psi(\Delta v, \Delta v) \geq 0$.  Since $|\Psi(\Delta v, v)|
\leq \epsilon \, \opnorm{\Sampcov}$, we obtain the lower bound
\begin{align*}
\sigmin \left ( \Sampcov \right )
= \inf_{v \in \Sphere{\pdim}} v^T \Sampcov v & \;
\geq \; \min_{v^b \in \Coverb} \Psi(v^b,v^b) - 2 \epsilon \opnorm{
\sampcov }.
\end{align*}
By the previously established upper bound\eqref{EqnARSig}(a), have
$\opnorm{ \sampcov } \leq \frac{24 \opnorm{\Sigma}}{(1-\maxop)}$ with
high probability.  Hence, choosing $\epsilon = \frac{(1-\maxop)
\sigmin(\Sigma)}{200 \opnorm{\Sigma}}$ ensures that \mbox{$2 \epsilon
\opnorm{\sampcov} \leq \sigmin(\Sigma)/4$.}

Consequently, it suffices to lower bound the minimum over the covering
set.  We first establish a concentration result for the function
$\Psi(v,v)$ that holds for any fixed $v \in \Sphere{\pdim}$.
Note that we can write
\begin{align*}
  \Psi(v, v) & = \frac{1}{\numobs} \sum_{i=1}^\numobs
  (\inprod{v}{\SpecX_i})^2,
\end{align*}
As before, if we define the random vector $Y \in \real^\numobs$ with
elements $Y_i = \inprod{v}{\SpecX_i}$, then $Y \sim N(0, R)$ with
$\opnorm{R} \leq \frac{2 \opnorm{\Sigma}}{1-\maxop}$.  Moreover, we
have $\trace(R)/\numobs = v^T \Sigma v \geq \sigmin(\Sigma)$.
Consequently, applying Lemma~\ref{LemMeta} yields
\begin{align*}
  \mprob \biggr[ \, \frac{1}{\numobs} \, \| Y \|_2^2 \; < \;
  \sigmin(\Sigma) - \frac{8 t \opnorm{\Sigma}}{1-\maxop} \biggr] &
  \leq 2 \exp \big (- \numobs (t- 2/\sqrt{\numobs})^2/2 \big ) + 2
  \exp( -\frac{\numobs}{2}),
\end{align*}
Note that this bound holds for any fixed $v \in \Sphere{\pdim}$.
Setting $t^* = \frac{(1-\maxop) \, \sigmin(\Sigma)}{16 \opnorm{\Sigma}}$
and applying the union bound yields that
\begin{align*}
\mprob \big[ \min_{v^b \in \Coverb} \Psi(v^b, v^b) \, < \,
  \sigmin(\Sigma)/2 \big] & \leq \; \big ( \frac{4}{\epsilon}
  \big)^\pdim \; \biggr \{ 2 \exp \big (- \numobs (t^*-
  2/\sqrt{\numobs})^2/2 \big ) + 2 \exp( -\frac{\numobs}{2}) \biggr
  \},
\end{align*}
which vanishes as long as $\numobs > \frac{4 \log
(4/\epsilon)}{(t^*)^2} \pdim$.


\subsection{Proof of Lemma~\ref{LemAutoregressiveNoise}}
\label{AppLemAutoregressiveNoise}

Let $\Sphere{\pdim} = \{ u \in \real^\pdim \, \mid \, \|u\|_2 = 1 \}$
denote the Euclidean sphere in $\pdim$-dimensions, and for positive
scalars $a$ and $b$, define the random variable $\Psi(a,b) \mydefn
\sup_{u \in a \, \Sphere{\pdim}} \sup_{v \in b \, \Sphere{\pdim}}
\inprod{\SpecX v}{\SpecV u}$.  Note that our goal is to upper bound
$\Psi(1,1)$.  Let $\Covera = \{u^1, \ldots, u^\covnuma \}$ and
$\Coverb = \{ v^1, \ldots, v^\covnumb \}$ denote $1/4$ coverings of
$\Sphere{\pdim}$ and $\Sphere{\pdim}$, respectively.  Following the
same argument as in the proof of Lemma~\ref{LemMultivarNoise}, we
obtain the upper bound
\begin{equation}
\label{EqnInterClaimTwo}
  \Psi(1,1) \; \leq \; 4 \max_{u^a \in \Covera, v^b \in \Coverb}
\inprod{\SpecX v^b}{\SpecV u^a}
\end{equation}
We now apply the union bound to control the discrete maximum.  It is
known (e.g.,~\cite{LedTal91,Matousek}) that there exists a $1/4$
covering of $\Sphere{\pdim}$ with at most $8^\pdim$ elements.
Consequently, we have
\begin{align}
\label{EqnUnionTwo}
\mprob \big[|\psi(1,1)| \geq 4 \delta \, \numobs \big] & \leq 8^{2
\pdim} \max_{u^a, v^b} \; \; \mprob \big[ \frac{|\inprod{\SpecX
v^b}{\SpecV u^a}|}{\numobs} \geq \delta \big].
\end{align}
It remains to obtain a tail bound on the quantity 
$\mprob \big[
\frac{|\inprod{\SpecX v}{\SpecV u}|}{\numobs} \geq \delta \big]$,
for any fixed pair \mbox{$(u, v) \in \Covera \times \Coverb$.}

For each $i=1, \ldots, \numobs$, let $\SpecX_i$ and $\SpecV_i$ denote
the $i^{th}$ row of $\SpecX$ and $\SpecV$.  Following some simple
algebra, we have the decomposition $\frac{\inprod{\SpecX v}{\SpecV
u}}{\numobs} = \Sah_1 - \Sah_2 - \Sah_3$, where
\begin{align*}
\Sah_1 & = \frac{1}{2 n} \sum_{i=1}^n \big( \binprod{u}{\SpecV_i} +
  \binprod{v}{\SpecX_i} \big)^2 - \frac{1}{2} (\nvar^2 + v^T \Sigma v) \\
\Sah_2 & =  \frac{1}{2\numobs} \sum_{i=1}^\numobs \big (
  \binprod{u}{\SpecV_i} \big )^2 - \nvar^2/2 \\
\Sah_3 & = \frac{1}{2\numobs} \sum_{i=1}^\numobs \big
  (\binprod{v}{\SpecX_i} \big )^2 - \frac{1}{2} v^T \Sigma v
\end{align*}
We may now bound each $\Sah_j, j = 1,2,3$ in turn; in doing so, we
make repeated use of Lemma~\ref{LemMeta}, which provides concentration
bounds for a random variable of the form $\|Y\|_2^2$, where \mbox{$Y
\sim N(0, Q)$} for some matrix $Q \succeq 0$.

\noindent \paragraph{Bound on $\Sah_2$:} We begin with $\Sah_2$, which
the easiest to control since (up to scaling by $\nvar$), it
corresponds to the deviation away from the mean of $\chi^2$-variable
with $\numobs$ degrees of freedom.  Consequently, applying
Lemma~\ref{LemMeta} with $Q = I$, we obtain
\begin{align}
\label{EqnSahTwo}
\mprob\big[ | \Sah_2 | \, > \, 4 \nvar^2 \,t \big] & \leq 2 \exp \Big
 (-\frac{\numobs \, (t - 2/\sqrt{\numobs})^2}{2} \Big ) + 2 \exp(
 -\numobs/2).
\end{align}

\paragraph{Bound on $\Sah_3$:} 
We can write the term $\Sah_3$ as a deviation of $\|Y\|_2^2/\numobs$
from its mean, where in this case the covariance matrix $Q$ is no
longer the identity.  In concrete terms, let us define a random vector
$Y \in \real^\numobs$ with elements $Y_i = \inprod{v}{\SpecX_i}$.  As
seen in the proof of Lemma~\ref{LemAutoregressiveRSC} from
Appendix~\ref{AppLemAutoregressiveRSC}, the vector $Y$ is zero-mean
Gaussian with covariance matrix $R$ such that $\opnorm{R} \leq \frac{2
\opnorm{\Sigma}}{1-\maxop}$ (see equation~\eqref{EqnRBound}).  Since
we have $\trace(R)/\numobs = v^T R v$, applying Lemma~\ref{LemMeta}
yields that
\begin{align}
\label{EqnSahThree}
\mprob \big[ | \Sah_3 | \, \geq \, \frac{8 \opnorm{\Sigma}}{1-\maxop}
t \big] & \leq \; 2 \exp \Big ( - \frac{\numobs \, (t \, -\,
2/\sqrt{\numobs})^2}{2} \Big ) + 2 \exp( -\numobs/2).
\end{align}

\paragraph{Bound on $\Sah_1$:}  To control this quantity, let us
 define a zero-mean Gaussian random vector $Z \in \real^\numobs$ with
elements $Z_i \, = \, \inprod{v}{\SpecX_i} \, + \,
\inprod{u}{\SpecV_i}$.  This random vector has covariance matrix
$S$ with elements
\begin{equation*}
S_{ij} = \Exs[Z_i Z_j] = \nvar^2 \delta_{ij} + (1-\delta_{ij})
  \nvar^2 v^T (\thetastar)^{|i-j|-1} u + v^T (\thetastar)^{|i-j|}
  \Sigma v,
\end{equation*}
where $\delta_{ij}$ is the Kronecker delta for the event $\{i = j\}$.
As before, by symmetry of $S$, we have $\opnorm{S} \leq \max_{i=1,
\ldots, \numobs} \sum_{j=1}^\numobs |S_{ij}|$, and hence
\begin{align*}
\opnorm{S} & \leq \nvar^2 + \opnorm{\Sigma} + \sum_{j=1}^{i-1} |\nvar^2
v^T (\thetastar)^{|i-j|-1} u + v^T (\thetastar)^{|i-j|} \Sigma v| +
\sum_{j=i+1}^n |\nvar^2 v^T (\thetastar)^{|i-j|-1} u + v^T
(\thetastar)^{|i-j|} \Sigma v| \\
& \leq \nvar^2 + \opnorm{\Sigma} + 2 \sum_{j=1}^{\infty} \nvar^2
  r^{j-1} + 2 \sum_{j=1}^\infty \matsnorm{\Sigma}{2} r^j \\
 & \leq \nvar^2 + \opnorm{\Sigma} + \frac{2 \nvar^2}{1-\maxop} +
  \frac{2 \maxop \opnorm{\Sigma}}{1-\maxop}.
\end{align*}
Morever, we have $\trace(S)/\numobs = \nvar^2 + v^T \ThetaStar v$, so
that by applying Lemma~\ref{LemMeta}, we conclude that
\begin{equation}
\label{EqnSahOne}
  \mprob \biggr[|\Sah_1 | \, > \, \big(\frac{12 \nvar^2}{1-\maxop} +
  \frac{12 \opnorm{\Sigma}}{1-\maxop} \big)
 \, t \biggr] \; \leq \;
 2
  \exp \Big ( - \frac{\numobs \, (t \, -\, 2/\sqrt{\numobs})^2}{2}
  \Big) + 2 \exp(-\numobs/2),
\end{equation}
which completes the analysis of this term.

Combining the bounds~\eqref{EqnSahThree},~\eqref{EqnSahTwo}
and~\eqref{EqnSahOne}, we conclude that for all $t > 0$,
\begin{align}
\label{EqnFixedUVBound}
\mprob \big[ \frac{|\inprod{\SpecX v}{\SpecV u}|}{\numobs} \geq
\frac{20 (\opnorm{\Sigma} + \nvar^2) \, t}{1-\maxop} \big] & \leq 6
\exp \Big ( - \frac{\numobs \, (t \, -\, 2/\sqrt{\numobs})^2}{2} \Big)
+ 6 \exp(-\numobs/2).
\end{align}
Setting $t = 10 \sqrt{\pdim/\numobs}$ and combining with the
bound~\eqref{EqnUnionTwo}, we conclude that
\begin{align*}
\mprob \big[|\psi(1,1)| \geq \frac{400 (\opnorm{\Sigma} +
\nvar^2)}{1-\maxop} \sqrt{\frac{\pdim}{\numobs}} \big] & \leq
8^{2 \pdim} \: \big \{ 6 \exp(-16 \pdim) +
6 \exp(-\numobs/2) \big \} \; \leq 12 \exp(-\pdim) 
\end{align*}
as long as $\numobs > ((4 \log 8) +1) \pdim$.


\section{Proof of Proposition~\ref{LemCompressedRSC}}
\label{AppLemCompressedRSC}

Note that $\|\Xop(\Theta)\|_2 = \sup_{u \in \Sphere{\bigN}}
\inprod{u}{\Xop(\Theta)}$, and that since the
claim~\eqref{EqnCompressedRSC} is invariant to rescaling, it suffices
to prove it for all $\Theta \in \real^{\lra \times \lrb}$ with
$\matsnorm{\Theta}{F} = 1$.  Letting $\arbrad \geq 1$ be a given
radius, we seek lower bounds on the quantity
\begin{equation*}
Z^*(\arbrad) \mydefn \inf_{\Theta \in \SpecSet(\arbrad)} \sup_{u \in
  \Sphere{\bigN}} \inprod{u}{\Xop(\Theta)}, \qquad \mbox{where
  $\SpecSet(\arbrad) = \{ \Theta \in \real^{\lra \times \lrb} \, \mid
  \, \matsnorm{\Theta}{F} = 1, \; \nucnorm{\Theta} \leq \arbrad \}$.}
\end{equation*}
In particular, our goal is to prove that for any $\arbrad \geq 1$, the
lower bound
\begin{equation}
\label{EqnInterBound}
\frac{Z^*(\arbrad)}{\sqrt{\bigN}} \geq \conone - \contwo \big[
  \frac{\lra + \lrb}{\bigN} \big]^{1/2} \;\arbrad
\end{equation}
holds with probability at least $1 - c_1 \exp(-c_2 \bigN)$.  By a
standard peeling argument (see Raskutti et al.~\cite{RasWaiYu09} for
details), this lower bound implies the claim~\eqref{EqnCompressedRSC}.

We establish the lower bound~\eqref{EqnInterBound} using Gaussian
comparison inequalities~\cite{LedTal91} and concentration of measure
(see Lemma~\ref{LemConcentration}).  For each pair $(u, \Theta) \in
\Sphere{\bigN} \times \SpecSet(\arbrad)$, consider the random
variable $Z_{u, \Theta} = \inprod{u}{\Xop(\Theta)}$, and note that it
is Gaussian with zero mean.  For any two pairs $(u, \Theta)$ and $(u',
\Theta)$, some calculation yields
\begin{align}
\label{EqnZCov}
\Exs \big[ (Z_{u,\Theta} - Z_{u',\Theta'})^2 ] & = 
\matsnorm{u \otimes
  \Theta - u' \otimes \Theta'}{F}^2.
\end{align}
We now define a second Gaussian process $\{Y_{u, \Theta} \, \mid \,
(u, \Theta) \in \Sphere{\bigN} \times \SpecSet(\arbrad) \}$ via
\begin{equation*}
Y_{u, \Theta} \mydefn \inprod{g}{u} + \tracer{G}{\Theta},
\end{equation*}
where $g \in \real^\bigN$ and $G \in \real^{\lra \times \lrb}$ are
independent with i.i.d. $N(0,1)$ entries.  By construction, $Y_{u,
  \Theta}$ is zero-mean, and moreover, for any two pairs $(u, \Theta)$
and $(u', \Theta')$, we have
\begin{align}
\label{EqnYcov}
\Exs \big[ (Y_{u,\Theta} - Y_{u',\Theta'})^2 ] & = 
\|u - u'\|_2^2 +
\matsnorm{\Theta - \Theta'}{F}^2.
\end{align}
It can be shown that for all pairs $(u, \Theta), (u', \Theta') \in
\Sphere{\bigN} \times \SpecSet(\arbrad)$, we have
\begin{align}
\matsnorm{u \otimes \Theta - u' \otimes \Theta'}{F}^2 & \leq \|u -
u'\|_2^2 + \matsnorm{\Theta - \Theta'}{F}^2.
\end{align}
Moreover, equality holds whenever $\Theta = \Theta'$.  The conditions
of the Gordon-Slepian inequality~\cite{LedTal91} are satisfied, so
that we are guaranteed that
\begin{align}
\Exs[\inf_{\Theta \in \SpecSet(\arbrad)} \|\Xop(\Theta)\|_2 ]\; = \;
\Exs \biggr[\inf_{\Theta \in \SpecSet(\arbrad)} \sup_{u \in
    \Sphere{\numobs}} Z_{u, \Theta} \biggr] & \geq \Exs
\biggr[\inf_{\Theta \in \SpecSet(\arbrad)} \sup_{u \in
    \Sphere{\numobs}} Y_{u, \Theta} \biggr]
\end{align}
We compute
\begin{align*}
\Exs \biggr[\inf_{\Theta \in \SpecSet(\arbrad)} \sup_{u \in
    \Sphere{\bigN}} Y_{u, \Theta} \biggr] & = \Exs \biggr [ \sup_{u \in
    \Sphere{\bigN}} \inprod{g}{u} \biggr] + \Exs \biggr [\inf_{\Theta \in
    \SpecSet(\arbrad)} \tracer{G}{\Theta} \biggr ] \\
& = \Exs[\|g\|_2] - \Exs[\sup_{\Theta \in \SpecSet(\arbrad)}
  \tracer{G}{\Theta}] \\
& \geq \frac{1}{2} \sqrt{\bigN} - \arbrad \, \Exs[\opnorm{G}].
\end{align*}
Since $G \in \real^{\lra \times \lrb}$ has i.i.d. $N(0,1)$ entries,
standard random matrix theory~\cite{DavSza01} implies that
$\Exs[ \, \opnorm{G} \,] \leq \sqrt{\lra} + \sqrt{\lrb}$.  Putting together
the pieces, we conclude that
\begin{align*}
\Exs[\inf_{\Theta \in \SpecSet(\arbrad)}
  \frac{\|\Xop(\Theta)\|_2}{\sqrt{\bigN}} ] & \geq \frac{1}{2} -
\frac{\sqrt{\lra} + \sqrt{\lrb}}{\sqrt{\bigN}} \, \arbrad.
\end{align*}

Finally, we need to establish sharp concentration around the mean.
Note that the function $f(\Xop) \mydefn \inf_{\Theta \in
\SpecSet(\arbrad)} \|\Xop(\Theta)\|_2/\sqrt{\bigN}$ is Lipschitz
with constant $1/\sqrt{\bigN}$, so that Lemma~\ref{LemConcentration}
implies that
\begin{align*}
\mprob \biggr[\inf_{\Theta \in \SpecSet(\arbrad)}
  \frac{\|\Xop(\Theta)\|_2}{\sqrt{\bigN}} \leq \frac{1}{2} - \arbrad
  \, \frac{\sqrt{\lra} + \sqrt{\lrb}}{\sqrt{\bigN}} - \delta \biggr] &
  \leq 2 \exp(-\bigN \delta^2/2) \qquad \mbox{for all $\delta > 0$.}
\end{align*}
Setting $\delta = 1/4$ yields the claim.


\section{Some useful concentration results}

The following lemma is classical~\cite{LedTal91,Massart03}, and yields
sharp concentration of a Lipschitz function of Gaussian random
variables around its mean.
\blems
\label{LemConcentration}
Let $X \in \real^{\numobs}$ have i.i.d. $N(0,1)$ entries, and let and
$f : \real^\numobs \to \real$ be Lipschitz with constant $L$ (i.e.,
$|f(x)-f(y)| \, \leq \, L \| x - y\|_2$ $\forall x, y \in
\real^\numobs$).  Then for all $t > 0$, we have
\begin{equation*}
  \mprob [ |f(X) - \ex f(X) | \, > \, t] \; \leq
  \; 2 \exp{\big (-\frac{t^2}{2 L^2} \big )}.
\end{equation*}
\elems

\noindent By exploiting this lemma, we can prove the following result,
which yields concentration of the squared $\ell_2$-norm of an
arbitrary Gaussian vector: \\

\blems
\label{LemMeta}
Given a Gaussian random vector $Y \sim N(0, Q)$, for all $t >
2/\sqrt{\numobs}$, we have
\begin{equation}
  \mprob \biggr[ \frac{1}{\numobs} \big | \| Y \|_2^2 - \trace{Q} \big
  | > 4 \, t \, \opnorm{Q} \biggr] \; \leq \; 
  2 \exp{\left (-
  \frac{\numobs (t- \frac{2}{\sqrt{\numobs}})^2}{2} \right)} + 2
  \exp{(-\numobs/2)}.
\end{equation}
\elems
\begin{proof}
Let $\sqrt{Q}$ be the symmetrix matrix square root, and consider the
function \mbox{$f(x) = \|\sqrt{Q} x\|_2/\sqrt{\numobs}$.}  Since it is
Lipschitz with constant $\opnorm{\sqrt{Q}}/\sqrt{\numobs}$,
Lemma~\ref{LemConcentration} implies that
\begin{equation}
  \label{EqnLipQ}
  \mprob \big[ \big | \; \|\sqrt{Q} X\|_2 - \ex \|\sqrt{Q} X\|_2 \;
  \big | > \sqrt{\numobs} \delta \big ] \leq 2 \exp \left (-\frac{\numobs
  \delta^2}{2 \opnorm{Q}} \right ) \qquad \mbox{for all $\delta > 0$.}
\end{equation}
By integrating this tail bound, we find that the variable $Z =
\|\sqrt{Q} X\|_2/\sqrt{\numobs}$ satisfies the bound $\var(Z) \leq 4
\opnorm{Q}/\numobs$, and hence conclude that
\begin{align}
\label{EqnVarQBound}
\big| \sqrt{\Exs[Z^2]} - |\Exs[Z]| \big| & = \big|
\sqrt{\trace(Q)/\numobs} - \Exs[\|\sqrt{Q} X\|_2/\sqrt{\numobs}]\big|
\; \leq \; \frac{2 \sqrt{\opnorm{Q}}}{\sqrt{\numobs}}.
\end{align}
Combining this bound with the tail bound~\eqref{EqnLipQ}, we conclude
that
\begin{equation}
  \label{EqnLipQTwo}
  \mprob \Big [ \frac{1}{\sqrt{\numobs}} \big | \|\sqrt{Q} X\|_2 -
  \sqrt{\trace(Q)} \; \big | > \delta + 2
  \sqrt{\frac{\opnorm{Q}}{\numobs}} \: \Big ] \leq 2
  \exp \left (-\frac{\numobs \delta^2}{2 \opnorm{Q}} \right )
  \qquad \mbox{for all
  $\delta > 0$.}
\end{equation}

Setting $\delta = (t-2/\sqrt{\numobs}) \, \sqrt{\opnorm{Q}} $ in the
bound~\eqref{EqnLipQTwo} yields that
\begin{align}
\label{EqnNegTerm}
  \mprob \Big [ \frac{1}{\sqrt{\numobs}} \big | \|\sqrt{Q} X\|_2 -
  \sqrt{\trace(Q)} \; \big | > t \, \sqrt{\opnorm{Q}} \: \Big ] \leq 2
  \exp \left (-\frac{\numobs (t-2/\sqrt{\numobs})^2}{2} \right ).
\end{align}
Similarly, setting $\delta = \sqrt{\opnorm{Q}}$ in the
tail bound~\eqref{EqnLipQTwo} yields that with probability greater
than $1 - 2 \exp(-\numobs/2)$, we have
\begin{align}
\label{EqnPlusTerm}
\biggr | \frac{\|Y\|_2}{\sqrt{\numobs}} +
\sqrt{\frac{\trace(Q)}{\numobs}} \biggr| & \leq
\sqrt{\frac{\trace(Q)}{\numobs}} + 3 \sqrt{\opnorm{Q}}
\; \leq \; 4 \sqrt{\opnorm{Q}}.
\end{align}
Using these two bounds, we obtain
\begin{align*}
\biggr | \frac{\|Y\|_2^2}{\numobs} - \frac{\trace(Q)}{\numobs} \biggr|
& = \biggr | \frac{\|Y\|_2}{\sqrt{\numobs}} -
\sqrt{\frac{\trace(Q)}{\numobs}} \biggr| \; \biggr |
\frac{\|Y\|_2}{\sqrt{\numobs}} + \sqrt{\frac{\trace(Q)}{\numobs}}
\biggr| \; \leq \; 4 \, t \, \opnorm{Q}
\end{align*}
with the claimed probability.
\end{proof}






\end{document}